\documentclass[final,leqno,onefignum,onetabnum]{siamltex1213}

\usepackage{bm,amssymb,amsmath,enumerate,color,txfonts,mathrsfs,tipa,multirow,epsfig,indentfirst}%
\usepackage{algorithm,float}%amsfonts,

\def\be{\begin{equation}}
\def\en{\end{equation}}
\def\beq{\begin{eqnarray}}
\def\eq{\end{eqnarray}}
\def\beqx{\begin{eqnarray*}}
\def\eqx{\end{eqnarray*}}

\makeatletter
  \newcommand\figcaption{\def\@captype{figure}\caption}
  \newcommand\tabcaption{\def\@captype{table}\caption}
\makeatletter
\newtheorem{remark}{Remark}[section]

\title{A discontinuous plane wave neural network method for Helmholtz equation and time-harmonic Maxwell's equations% with large wave numbers
 \thanks{The first author was supported by Shandong Provincial Natural Science Foundation under the grant ZR2020MA046.
The second author was supported by the Natural Science Foundation of China G12071469.}}
\author{LONG YUAN
\thanks{College of Mathematics and Systems Science, Shandong University of Science and Technology,
Qingdao 266590, China (sdbjbjsd@163.com).}
 \and QIYA Hu
\thanks{Corresponding author. 1. LSEC, Institute of Computational Mathematics and
Scientic/Engineering Computing, Academy of Mathematics and Systems
Science, Chinese Academy of Sciences, Beijing 100190, China; 2. School of Mathematical Sciences, University of Chinese Academy
of Sciences, Beijing 100049, China
(hqy@lsec.cc.ac.cn).}
 }

\begin{document}

\maketitle

\begin{abstract} %with large wave numbers
In this paper we propose a {\it discontinuous} plane wave neural network (DPWNN) method  with $hp-$refinement for approximately solving Helmholtz equation and time-harmonic Maxwell equations. In this method, we define a quadratic functional as in the plane wave least square (PWLS) method with $h-$refinement and introduce new discretization sets spanned by element-wise neural network functions with a single hidden layer, where the activation function on each element is chosen as a complex-valued exponential function like the plane wave function. The desired approximate solution is recursively generated by iteratively solving a quasi-minimization problem associated with the functional and the sets described above, which is defined by a sequence of
%approximations adaptively corrected by searching
approximate minimizers of the underlying residual functionals, where plane wave direction angles and activation coefficients are alternatively computed by iterative algorithms. For the proposed DPWNN method, the plane wave directions are adaptively determined in the iterative process, which is different from that in the standard PWLS method (where the plane wave directions are preliminarily given). Numerical experiments will confirm that this DPWNN method can generate approximate solutions with higher accuracy than the PWLS method.
\end{abstract}

\begin{keywords}
Helmholtz equations, time-harmonic Maxwell's equations, least squares, element-wise neural network, plane wave activation function, iterative algorithms, convergence
\end{keywords}
\begin{AMS}
65N30, 65N55.
\end{AMS}

\pagestyle{myheadings} \thispagestyle{plain} \markboth{LONG YUAN AND QIYA HU}{Discontinuous plane wave neural networks for Helmholtz equation and Maxwell's equations}

\section{Introduction}

Acoustic, elastic and electromagnetic propagation problems arise in many areas of physical and engineering interest, in areas as diverse as radar, sonar, building acoustics, medical and seismic imaging. Mathematically, the problem of wave propagation is often modeled by Helmholtz equation or Maxwell's equations (see \cite{ALZOU, BZZOU, Imbert, zhao}). One of the main difficulties in the numerical simulation of the propagating problem arises from the usually oscillatory nature of the solutions. Plane wave methods have been designed for finite element discretization of time-harmonic wave equations (for example, the Helmholtz and Maxwell's equations) in the last years, two typical plane wave methods are the plane wave discontinuous Galerkin (PWDG) method \cite{Git, ref21,pwdg} and the plane wave least squares (PWLS) method \cite{hy2,hy3,ref12,peng}. The plane wave methods have an important advantage over Lagrange finite elements for discretization of the Helmholtz equation and time-harmonic Maxwell equations: to achieve the same accuracy, relatively smaller number of degrees of freedom are enough in the plane wave-type methods owing to the particular choice of the basis functions that satisfy the considered partial differential equations (PDEs) without boundary conditions. However the standard plane wave methods are heavily
unstable \cite{ref11,ref21,HMPsur,HGA}, which limits real accuracies of the resulting approximate solutions.

%Recent works \cite{BGL, yuan2} also have led to the development of algorithms that incorporate the highly oscillatory nature of the solutions into the basis functions.

The neural network method has been widely applied to approximately solving PDEs. To obtain an approximate solution of the considered PDE via neural networks, the primary key step is to choose suitable neural networks for minimizing the PDE residual, and several approaches have been proposed to accomplish this \cite{berg, eyu, hexu, Kharazmi, raissi, zang}. While such physics-informed neural networks and variational approaches have enjoyed success in particular cases, in practice, the achievable relative error is prone to stagnate around $0.1-1.0\%$ no matter how many neurons or layers are used to define the underlying network architecture.
To overcome this plight, {\it continuous} Galerkin neural networks, based on the adaptive construction of a sequence of finite-dimensional subspaces whose basis functions are realizations of a sequence of neural networks, have been explored for variational problems with symmetric, bounded and coercive bilinear forms in \cite{AD}.
%This approach relies on the existing assumption that the bilinear form of the associated variational problem satisfies {\it symmetric, bounded} and {\it coercive} properties.
The sequential nature of the algorithm offers a systematic approach to enhancing the accuracy of a given approximation, and provides a useful indicator for the error in the energy norm that can be used as a criterion for terminating the sequential updates. The convergence results established in \cite{AD} rely on a key assumption that the neural network parameters are bounded.

There are a few works to apply the neural network method to the discretization of the Helmholtz equation, for example, the work \cite{chung} developed a deep learning approach to learn ray directions at discrete locations by analyzing highly oscillatory wave fields, and used the plane wave basis functions with the approximate ray directions in the (interior-penalty) discontinuous Galerkin plane wave method to solve the corresponding Helmholtz equations at higher frequencies. In the present paper, we try to combine the neural network method with the PWLS method for the discretization of Helmholtz equation and Maxwell equations.
We define a sequence of sets of {\it discontinuous} plane wave neural networks, where the activation function on each element is chosen as $e^{i\omega x}$ ($\omega$ denotes the wave number) and a single hidden layer is considered, and apply it to the discretization of a minimization problem as in the PWLS method for Helmholtz equation and time-harmonic Maxwell equations. The desired approximate solution is recursively generated by searching approximate {\it quasi-minimizers} on the sequence of sets for the residual functionals defined by the former {\it quasi-minimizer}, where plane wave direction angles and activation coefficients are alternatively computed by iterative algorithms. For convenience, we call the proposed method as DPWNN method.
%The DPWNN method described above falls into the class of Galerkin neural networks, but with new activation function.
The new choice of activation function can enhance the accuracies of the resulting approximate solutions for the underlying models and accelerate the convergence of the iterative algorithms.
We establish convergence results of the DPWNN method without the assumption on the boundedness of the neural network parameters. The performance of the algorithms will confirm the effectiveness of the proposed method.
The new DPWNN method differs from the Galerkin neural network method introduced in \cite{AD} both in algorithm design and theoretical analysis.

%\textcolor{red}{
%The new DPWNN method differs from the existing Galerkin neural network introduced in \cite{AD} in the sense that: i) the strictly constrained variational formulation is replaced by the adaptively constructed quadratic functional %and a relaxed notion of quasi-minimization problem is introduced;  ii) there is no constraint of the boundedness imposed on the neural network parameters; iii) the DPWNN method is simpler and easier to implement. The desired %approximate solution for our DPWNN method is recursively generated by directly searching approximate minimizers on the sequence of sets for the residual functionals defined by the former minimizer, while the work \cite{AD} used %continuous neural networks to approximate variational equations, based on the adaptive construction of a sequence of finite-dimensional subspaces whose basis functions are realizations of a sequence of neural networks by %learning the dual representation of the weak residual.
%}

 %new discretization sets spanned by element-wise neural network functions are introduced in the new objective functional so that by choosing $h\approx \mathcal{O}(\frac{4\pi}{\omega})$ and gradually increasing the width $n_r$ of the networks,  the performance of the iterative algorithms is mildly dependent on wave numbers.

Differently from the standard PWLS method, for the new DPWNN method the plane wave directions are not fixed and can be corrected in the iterative process by searching approximate minimizers of the involved residuals, and the approximate solutions are adaptively computed by augmenting the discontinuous plane wave neural networks (to gradually decrease the residuals). Numerical results in section 6 will confirm that the DPWNN method can generate
approximate solutions with higher accuracy than the PWLS method.

%Besides, compared with the deep neural networks for the wave equations \cite{Kharazmi, Karniadakis}, the DPWNN method is more intuitive and easier to implement, yet two different types of methods both can
%achieve high accuracies.

%\textcolor{red}{}

The paper is organized as follows: In section 2, we define a minimization problem for the underlying Helmholtz equation. In section 3, we introduce discontinuous plane wave neural networks. In section 4, we describe
a discretization for the minimization problem, and design three iterative algorithms to generate approximate solutions of the discrete minimization problem, and give an error estimate of the approximate solutions.
In section 5, we apply the proposed DPWNN method for the discretization of the time-harmonic Maxwell equations. Finally, we report some numerical results to confirm the
effectiveness of the proposed method in section 6.

\section{A minimization method for Helmholtz equations}
Let $\Omega \subset \mathbb{R}^d~(d=2, 3)$ be a bounded polygonal/polyhedral domain with boundary $\gamma=\partial\Omega$. Consider Helmholtz equations which is formalized, normalizing the wave's velocity to 1, by
\begin{eqnarray}
\left\{\begin{array}{ll} -\Delta u-\omega^2u=0& \text{in}\quad
\Omega,\\
(\partial_\text{\bf n} + \text{i}\omega)u=g& \text{on}\quad\gamma=\partial\Omega,
\end{array}\right.
\label{helm1}
\end{eqnarray}
where $g \in L^2(\partial\Omega)$. The outer normal derivative is referred to by $\partial_{\bf n}$ and the angular frequency by $\omega$.

%\textcolor{red}{single subdomain in numerical tests}

The considered variational formulation is based on a triangulation of the solution domain. Let \(\Omega\) be divided into a partition
in the sense that
$$ \overline{\Omega}=\bigcup_{k=1}^N\overline{\Omega}_k,\quad
\Omega_k\bigcap\Omega_j=\emptyset
 \quad\text{ for }k\not=j. $$
Let \( {\cal T}_h\) denote the triangulation comprising the elements
\(\{\Omega_k\}\), where \(h\) is the meshwidth of the triangulation.
Define $$ \Gamma_{kj}=\partial\Omega_k\bigcap\partial\Omega_j
\quad\text{for
 }k\not=j$$
 and
$$ \gamma_k=\overline{\Omega}_k\bigcap\partial\Omega
 \quad(k=1,\ldots,N),~~
 \gamma=\bigcup^N_{k=1}\gamma_k.$$

Set $u|_{\Omega_k}=u_k$ $(k=1,\cdots,N$). Then the reference problem
to be solved consists in finding local acoustic pressures
$u_k\in H^1(\Omega_k)$ such that

\begin{eqnarray}
\left\{\begin{array}{ll} -\Delta u_k-\omega^2u_k=0& \text{in}\quad
\Omega_k,\\
(\partial_\text{\bf n} + \text{i}\omega)u_k = g
&\text{on}\quad\gamma_k,
\end{array}\right.\quad\quad(k=1,2,\ldots,N),
\label{helm2}
\end{eqnarray}
and
\begin{eqnarray}
\left\{\begin{array}{ll}u_k-u_j=0 & \text{over}\quad \Gamma_{kj},\\
\iffalse {\partial u_k\over \partial \text{{\bf n}}_k}+{\partial
u_j\over\partial {\bf n}_j}=0& \text{over}\quad \Gamma_{kj},\\\fi
\partial_{\text{\bf n}_k}u_k+\partial_{\text{\bf n}_j}u_j=0& \text{over}\quad
\Gamma_{kj}.
\end{array}\right.\quad\quad(k\neq j;~ k,j=1,2,\cdots,N).
\label{helm3}
\end{eqnarray}

Let $V(\Omega_k)$ denote the space of the functions which verify the homogeneous Helmholtz equation on each element $\Omega_k$:
\be
 V(\Omega_k) = \bigg\{v_k \in H^1(\Omega_k); -\Delta v_k - \omega^2 v_k = 0 \bigg\}.
 \en
Define
\be
V({\cal T}_h) = \{ v\in L^2(\Omega): ~v|_{\Omega_k}\in V(\Omega_k) ~\forall \Omega_k\in {\cal T}_h  \}.
\en
%\section{The PWLS method}
%Set $v|_{\Omega_k}=v_k$ $(k=1,\cdots,N$). Define a functional

Based on the above observation, we define a functional
\begin{eqnarray} \label{funct}
 J(v) &=& \sum_{k=1}^N\int _{\gamma_k} |(\partial_{\bf
n}+i\omega )v_k-g|^2 ds  \cr
& + & \sum_{j\not=k} \bigg( \alpha \int_{\Gamma_{kj}}
|v_k-v_j|^2 ds + \beta \int_{\Gamma_{kj}} | (\partial_{\text{\bf n}_k}
v_k+\partial_{\text{\bf n}_j} v_j ) |^2ds \bigg), ~~\forall v\in V({\mathcal T}_h),
\end{eqnarray}
where $\alpha$ and $\beta$ are given positive numbers. The introduction of the two relaxation parameters $\alpha$ and $\beta$ aims to balance the term $\alpha \int_{\Gamma_{kj}} |v_k-v_j|^2 ds$ and the term $\beta \int_{\Gamma_{kj}} | (\partial_{\text{\bf n}_k} v_k+\partial_{\text{\bf n}_j} v_j ) |^2ds$ in the functional $J$, especially when the wave number $\omega$ is large, which results in the analytic solution $u$ becomes high oscillating. In general we choose $\alpha=\omega^2$ and $\beta=1$ (see the detailed in \cite{hy}).

Consider the minimization problem: find $u\in {V}({\cal T}_h)$ such that
\begin{equation} \label{minmum}
J({u})=\min\limits_{{v}\in {V}({\cal T}_h)}~J(v).
\end{equation}

The variational problem associated with the minimization problem (\ref{minmum}) can be expressed as follows (refer to \cite{hy}): find $u \in V({\mathcal T}_h)$ such that
\begin{eqnarray}
 &&\sum_{k=1}^N\int _{\gamma_k}\big((\partial_{\bf
n}+i\omega )u_k-g\big)\cdot\overline{(\partial_{\bf n}+i\omega) v_k}ds+
 \sum_{j\not=k} \bigg(\alpha\int_{\Gamma_{kj}}
(u_k-u_j)\cdot\overline{(v_k-v_j)}ds \cr
&&\quad\quad+\beta\int_{\Gamma_{kj}}(\partial_{\text{\bf n}_k}
u_k+\partial_{\text{\bf n}_j} u_j
)\cdot\overline{(\partial_{\text{\bf n}_k} v_k+\partial_{\text{\bf
n}_j} v_j )} ds\bigg)=0, ~~\forall v\in V({\mathcal T}_h).
\label{va1}
\end{eqnarray}
Define the sesquilinear form $a(\cdot, \cdot)$ by
\beq
& a(u,v) = \sum_{k=1}^N\int _{\gamma_k} (\partial_{\bf
n}+i\omega )u_k \cdot\overline{(\partial_{\bf n}+i\omega) v_k}ds
+ \sum_{j\not=k} \bigg(\alpha\int_{\Gamma_{kj}}
(u_k-u_j)\cdot\overline{(v_k-v_j)} ds  \cr
& + \beta\int_{\Gamma_{kj}}(\partial_{\text{\bf n}_k}
u_k +\partial_{\text{\bf n}_j} u_j ) \cdot \overline{(\partial_{\text{\bf n}_k} v_k+\partial_{\text{\bf n}_j} v_j )} ds \bigg), ~\forall v\in V({\mathcal T}_h),
\eq
and semilinear form $L(v)$ by
\be
L(v)=\sum_{k=1}^N\int _{\gamma_k} g ~\overline{(\partial_{\bf
n}+i\omega )v_k} ds, ~\forall v\in V({\mathcal T}_h).
\en
 Then (\ref{va1}) can be written as the following variational formulation:
 \begin{eqnarray}\label{pwlsvar}
\left\{\begin{array}{ll} \text{Find} \ u\in V({\mathcal T}_h) ~~s.t. \\
a(u,v)=L(v) \quad \forall v\in V({\mathcal T}_h).
\end{array}\right.
\end{eqnarray}

 It is easy to see that $a(v,v)$ is a norm on $V({\mathcal T}_h)$ (see \cite{hy}).
For ease of notation, this norm is denoted by $||| v |||^2$. By the Cauchy-Schwarz inequality, we also have that $a(\cdot,\cdot)$ satisfies the boundedness: $|a(u,v)|\leq ||| u |||~||| v |||, ~\forall u,v \in V({\mathcal T}_h)$, and that $L(v)$ is bounded satisfying $|L(v)| \leq ||g||_{L^2(\gamma)} ~||| v |||$.  Furthermore, a standard application of the Risez representation theorem \cite{FO} implies the existence of a unique $u\in V({\mathcal T}_h)$ satisfying (\ref{pwlsvar}).

\section{Discontinuous plane wave neural networks}
A {\it discontinuous} neural network consists of a single hidden layer of $n\in \mathbb{R}$ neurons on each element $\Omega_k\in {\cal T}_h$. Define a function $\varphi^{\theta}: \mathbb{R}^d \rightarrow \mathbb{C}$ as follows:
\be \label{nn_app}
\varphi^{\theta}({\bf x})|_{\Omega_k} = \sum_{j=1}^n c_j^{(k)} \sigma(W_j^{(k)} \cdot {\bf x} + b_j^{(k)}), ~\forall \Omega_k \in {\cal T}_h,
\en
where $n$ is called the width of the network; $W_j^{(k)}\in \mathbb{R}^d, b_j^{(k)}\in \mathbb{R}$ are defined element-wisely; $c_j^{(k)} \in \mathbb{C}$ are element-wise coefficients; $\sigma: \mathbb{R}\rightarrow \mathbb{C}$ is a bounded activation function.
%For simplicity, we shall use the elementwise notation $D|_{\Omega_k}=[D_1, \cdots, D_n]\in \mathbb{R}^{d\times n}$ referred as weights, $b|_{\Omega_k}=[b_1, \cdots, b_n]\in \mathbb{R}^{n}$ referred as biases, and %$c|_{\Omega_k}=[c_1, \cdots, c_n]^T\in \mathbb{C}^n$. A set of nonlinear and linear parameters over a collection of all elements is denoted by $\theta=\{D,b,c\}$.

%\textcolor{red}{motivation provided to justify the need or benefit to use NNs (mention
%"neural network approximation" but without explaining the context).}

As pointed out in section 1, the plane wave method has obvious advantages over Lagrange finite elements for discretization of the time-harmonic wave equations. Thus, in this paper we choose a {\it plane wave type} function
$e^{i\omega x}$ as the activation function $\sigma(x)$ and set $W_j^{(k)}={\bf d}_j^{(k)}$, where ${\bf d}_j^{(k)}\in \mathbb{R}^d, |{\bf d}_j^{(k)}|=1$ are $n$ different propagation directions on
$k$-th element $\Omega_k\in {\cal T}_h$. Then we write
$$ \sigma(W_j^{(k)} \cdot {\bf x} + b_j^{(k)})=e^{i\omega({\bf d}_j^{(k)} \cdot {\bf x}+b_j^{(k)})}=e^{i\omega b_j^{(k)}}e^{i\omega{\bf d}_j^{(k)}\cdot {\bf x}},  \quad ({\bf d}_j^{(k)}\in \mathbb{R}^d, |{\bf d}_j^{(k)}|=1). $$
Notice that $e^{i\omega b_j^{(k)}}$ can be emerged into the coefficient $c_j^{(k)}$, we can ignore $b_j^{(k)}$ in the current situation and define the element-wise neural network
\be \label{nn_app2_general}
\varphi^{\theta}({\bf x})|_{\Omega_k} = \sum_{j=1}^n c_j^{(k)} \text{e}^{\text{i} \omega {\bf d}_j^{(k)} \cdot {\bf x}}, ~~{\bf x}\in \Omega_k\quad (|{\bf d}_j^{(k)}| = 1),
\en
with the element-wise training parameters $\theta^{(k)}=\{{\bf d}_j^{(k)}; c_j^{(k)}\}_{j=1}^n$ defined below. Let $V_n({\mathcal T}_h)$ be the set of all functions of the form (\ref{nn_app2_general}), namely,
\be \label{vsign_general}
V_n({\mathcal T}_h) := \bigg\{v\in L^2(\Omega): ~v|_{\Omega_k}=\sum_{j=1}^n c_j^{(k)} \text{e}^{\text{i} \omega {\bf d}_j^{(k)} \cdot {\bf x}}, ~~{\bf x}\in \Omega_k,~~\Omega_k\in {\mathcal T}_h~~(c_j^{(k)} \in \mathbb{C}, ~{\bf d}_j^{(k)}\in \mathbb{R}^d, ~|{\bf d}_j^{(k)}| = 1)\bigg\}.
\en
We emphasize that $V_n({\mathcal T}_h)$ is not a linear space since the directions $\{{\bf d}_j^{(k)}\}$ are not given. But, for simplicity of exposition,
we still call %the discontinuous Galerkin \textcolor{red}{subspace}
$V_n({\mathcal T}_h)$ as a plane wave neural network space.

\subsection{A universal approximation}

As stated above, we employ the plane wave basis functions as the activation functions, which generate the associated plane wave neural network. We first state the existing results from \cite[Corollary 5.5]{mhp}: the $hp-$approximation estimates for homogeneous Helmholtz solutions in $H^s(\Omega_k)$ with plane wave approximations. For ease of notation, we only give a simplified form of the result (see \cite[Lemma 4.4]{hy3}).

% (see \cite{mhp})
 Let \(s\) and \(m\) be given positive integers
satisfying \(m\geq 2(s-1)+1=2s-1\). Let the
number \(n\) of plane wave propagation directions be chosen as
\(n=2m+1\) in 2D and \(n=(m+1)^2\) in 3D, respectively.

\begin{lemma}  \label{elementappro}
Let $2\leq s\leq {m+1\over 2}$ with a sufficiently large $m$.
Let $u\in H^{s}(\Omega_k)$ be a solution of the homogeneous Helmholtz
equation for each element $\Omega_k$. Then, there exists $[c_1^{(k)},c_2^{(k)},\cdots,c_n^{(k)}] \in \mathbb{C}^n$ and $[{\bf d}_1^{(k)}, {\bf d}_2^{(k)},\cdots, {\bf d}_n^{(k)}]\in R^{d\times n}$,  such that
\be
||u - \sum_{l=1}^n c_l^{(k)} \text{e}^{\text{i} \omega {\bf d}_l^{(k)} \cdot{\bf x}} ||_{j,\Omega_k}\leq C h^{s-j}m^{-\lambda(s-j-\varepsilon)}||u||_{s,\omega,\Omega_k}~~~~(0\leq j\leq s; k=1,\cdots,N),\label{approx}
\en
where $\lambda>0$ is a constant depending only on the shape of the elements (in particular, $\lambda=1$ for the
 case of two dimensions), $\varepsilon=\varepsilon(m)>0$ satisfies $\varepsilon(m)\rightarrow 0$ when $m\rightarrow\infty$, and the $\omega-$weighted Sobolev norm is defined by
$||v||_{s,\omega,D} = \sum_{j=0}^{s}\omega^{2(s-j)}|v|_{j,D}^2.$
\end{lemma}

With the help of the above approximation property, we can further obtain the following universal approximation property of discontinuous plane wave neural networks, which can also be achieved by the classical approximate procedure \cite{HO}.

%\textcolor{red}{plane wave approximation as plane wave neural network}

 %\textcolor{red}{plane wave approximation on p}

\begin{theorem} \label{aninappro}
%\textcolor{red}{piecewise activation function, piecewise smoothness}
Suppose that $1\leq p< \infty, ~ 0\leq s < \infty$, and that $\Omega\subset \mathbb{R}^d$ is compact, % If $\sigma|_{\Omega_k}\in C^s(\Omega_k)$ is nonconstant and bounded,
 then $V_n({\mathcal T}_h)$ is dense in the Sobolev space
\be\nonumber
W^{s,p}({\cal T}_h) := \bigg\{v\in L^p(\Omega): D^{\alpha}v\in  L^p(\Omega_k)~~\forall |\alpha|\leq s, k=1,\cdots,N   \bigg\}.
\en
In particular, for any given function
$f \in  W^{s,p}({\cal T}_h)$ and $\tau  > 0$, there exists $n(\tau , f) \in  \mathbb{N}$  and $\tilde{f} \in  V_{n(\tau ,f )}({\cal T}_{h})$ such that $| | f-\tilde{f}| |_{W^{s,p}({\cal T}_h)} < \tau$.
\end{theorem}

%\textcolor{red}{define well-definedness}
\iffalse
 Although the discontinuous neural network $V_n^{h}$ is neither closed nor compact in $L^p(\Omega)$, by restricting the boundedness of network parameters (see \cite{PRV}), the discontinuous neural network with bounded parameters is both closed and compact in $L^p(\Omega)$.
 For the well-posedness of the discretized optimization problems (\ref{iteradisminmum}) below, the subset $V_{n,C}^{h}$ of $V_n^{h}$ consisting of functions with bounded parameters needs to be defined (refer to \cite{AD}).

 %\textcolor{red}{effect}

Given bounded constant $C>0$, we define
\be\label{bounnn}
V_{n,C}^{h}:= \bigg\{  v\in V_n^{h}({\cal T}_h): ||\theta||\leq C   \bigg\},
\en
where $||\theta||=\text{max}_k||\theta_k||$,  and $||\theta_k||:=\text{max}|\Phi_{j}|+\text{max}|c_j|$ for each element $\Omega_k\in {\cal T}_h$.
\fi

\begin{remark} Since the standard plane wave method is heavily unstable \cite{ref11,ref21,HMPsur,HGA}, the coordinate coefficients $\{c_l^{(k)}\}$ in (\ref{approx}) are unbounded when the number $n$ of basis functions
increases to infinity (see section 4 of \cite{Par2022}). This means that the universal approximation result Theorem \ref{aninappro} may not hold if the coefficients $\{c_j^{(k)}\}$ in the set $V_n({\mathcal T}_h)$
is required to be bounded. Because of this, we have not assumed the boundedness of the coefficients $\{c_j^{(k)}\}$ in the definition of the set $V_n({\mathcal T}_h)$. Then the set $V_n({\mathcal T}_h)$
is not closed (topologically) in ${V}({\cal T}_h)$, which make us have to consider a quasi-minimization problem (\ref{iteradisminmum}) instead of the standard minimization problem.
%and different analysis from that in \cite{AD}, where boundedness of neural network parameters was required.
\end{remark}

\subsection{Pre-definition of plane wave propagation directions}\label{helmcase}

For the 2D case, plane wave propagation directions ${\bf d}_{l}^{(k)}$ ($l=1,\cdots,n$) are defined by
\be \label{nn_app2}
{\bf d}_l^{(k)} = [cos(\vartheta_l^{(k)}), sin(\vartheta_l^{(k)})]^T, ~\vartheta_l^{(k)}\in (-\pi, \pi].
\en
A set of plane wave direction angles and activation coefficients defined on every element are denoted by $\vartheta|_{\Omega_k}=\vartheta^{(k)}=[\vartheta_1^{(k)}, \cdots, \vartheta_n^{(k)}]^T$,
 $c|_{\Omega_k}=c^{(k)}=[c_1^{(k)}, \cdots, c_n^{(k)}]^T\in \mathbb{C}^n$, respectively. A set of these parameters over a collection of all elements is denoted by $\theta=\{\vartheta; c\}$. In particular, a set of plane wave propagation angles defined on every element is collectively denoted by $\Phi=\{\vartheta\}$.

 For the 3D case, plane wave propagation directions ${\bf d}_{l}^{(k)}$ ($l=1,\cdots,n$) based on spherical coordinates are defined by
\begin{equation} \label{3ddirec}
{\bf d}_l^{(k)} = \left ( {\begin{array}{c}
sin(\zeta_m^{(k)}) cos(\vartheta_t^{(k)})  \\
sin(\zeta_m^{(k)}) sin(\vartheta_t^{(k)})  \\
cos(\zeta_m^{(k)})
\end{array}}
\right ),  ~\zeta_m^{(k)}\in [0, \pi],~m =1,\cdots, m^{\ast}; ~\vartheta_t^{(k)}\in (-\pi, \pi], ~t=1,\cdots, t^{\ast},
\end{equation}
where $m^{\ast}$ and $t^{\ast}$ are positive integers, which satisfy $t^{\ast}\approx 2m^{\ast}$ and $n=m^{\ast}t^{\ast}$. We will use the element-wise notation $\zeta|_{\Omega_k}=\zeta^{(k)}=[\zeta_1^{(k)}, \cdots, \zeta_{m^{\ast}}^{(k)}]^T\in \mathbb{R}^{m^{\ast}}$, $\vartheta|_{\Omega_k}=\vartheta^{(k)}=[\vartheta_1^{(k)}, \cdots, \vartheta_{t^{\ast}}^{(k)}]^T \in \mathbb{R}^{t^{\ast}}$, and $c|_{\Omega_k}=c^{(k)}=[c_1^{(k)}, \cdots, c_{n}^{(k)}]^T\in \mathbb{C}^{n}$. A set of these parameters over a collection of all elements is denoted by $\theta=\{\vartheta;\zeta; c\}$. In particular, a set of plane wave propagation angles defined on every element is collectively denoted by $\Phi=\{\vartheta; \zeta\}$.

Moreover, we would like to point out that $\zeta_m^{(k)}=0$ implies the resulting basis functions $\{e^{\text{i}\omega {\bf d}_{l}^{(k)}\cdot{\bf x}}\}$ are linearly dependent for different $\vartheta_t^{(k)}$ on each element $\Omega_k$. Similarly to the case of $\zeta_m^{(k)}=\pi$. To avoid this situation, on the beginning of training epoches in Algorithm 4.2, the angle $\zeta_m^{(k)}$ with the z-axis will not be initialized to 0 or $\pi$ in a cycle $T=\pi$. For each training epoch in Algorithm 4.2, the associated angles $\{\zeta_m^{(k)}\}$ satisfying $sin(\zeta_m^{(k)})\approx 0$ need to be corrected by adding suitable disturbance term $\varepsilon_{disturb}$ to guarantee the linear independence of the resulting basis functions $\{e^{\text{i}\omega {\bf d}_{l}^{(k)}\cdot{\bf x}}\}$ for different $\vartheta_t^{(k)}$.

%\textcolor{red}{for iterative case, we have the checking procedure.}

\iffalse
Subsequently, a neural network consists of a single hidden layer of $n\in \mathbb{R}$ neurons on each element $\Omega_q\in {\cal T}_h$, defining an output function $\varphi^{NN}^{NN}: \mathbb{R}^3 \rightarrow \mathbb{C}$ as follows:
\be \label{3dhelm_nn_app}
\varphi^{NN}^{NN}(x; \theta)|_{\Omega_q} = \sum_{m=1}^{n_1}\sum_{k=1}^{n_2}  c_l ~e^{\text{i}\omega
W_{l}\cdot x}, ~l=(m-1)n_2+k, % ~x\in \Omega_q \subset {\cal T}_h,
\en
\fi

%\input{chapter/chap-DPWNN}
\section{A discontinuous plane wave neural network iterative method}\label{helmdis}

In this section, we present iterative algorithms to generate discontinuous plane wave neural network solutions approximately satisfying the minimization problem (\ref{minmum}).
%inspired by the framework \cite{AD} of constructing the {\it continuous} Galerkin neural networks.
%we construct a {\it discontinuous} Galerkin plane wave neural networks for accelerating the convergence.
%, which naturally bring the $hp$-convergence of the proposed algorithm . %We would like to remark that

\subsection{Main results}

For the defined functional $J: { V}({\cal T}_h)\rightarrow \mathbb{R}$, we use $J'(v)\in (V({\cal T}_h))^{\ast}$ (dual space)
%: { V}({\cal T}_h)\rightarrow \mathbb{C}$
to represent the first-order Fr$\acute{e}$chet derivative of $J$ at $v\in { V}({\cal T}_h)$.
Similarly, we use $J''(v): V({\cal T}_h)\rightarrow (V({\cal T}_h))^{\ast}$
%{ V}({\cal T}_h)\times{ V}({\cal T}_h) \rightarrow \mathbb{C}$
to denote the second-order Fr$\acute{e}$chet derivative of $J$ at $v\in { V}({\cal T}_h)$.
Let $\langle \cdot, \cdot \rangle$ denote the duality pairing (which extends the standard $L^2$ inner product) between $V({\cal T}_h)$ and $(V({\cal T}_h))^{\ast}$.
It is easy to verify that
\be \nonumber
\langle J'(v), w \rangle = 2\text{Re}\big\{ a(v,w)-L(w) \big\},\quad v,w\in { V}({\cal T}_h)
\en
%\be \nonumber
%J'(v)w= a(v,w) - L(w),\quad v,w\in { V}({\cal T}_h)
%\en
and
\be \nonumber
\langle J''(v) w_1,w_2 \rangle = 2\text{Re} \big\{ a(w_1,w_2) \big\}, \quad v,w_1,w_2\in { V}({\cal T}_h).
\en
%\be \nonumber
%J''(v)(w_1,w_2) = a(w_1,w_2), \quad v,w_1,w_2\in { V}({\cal T}_h).
%\en

Suppose $u_{r-1}\in { V}({\cal T}_h)$ is an approximate solution of (\ref{minmum}). Let $\xi_r\in {V}({\cal T}_h)$ be the solution of
the minimization problem:
\begin{equation} \label{iteraminmum}
J(u_{r-1} + {\xi_r})=\min\limits_{{v}\in {V}({\cal T}_h)}~J(u_{r-1} + v).
\end{equation}
Then we have
\be \nonumber %\label{derizero}
\langle J'(u_{r-1} + \xi_r), v \rangle = 0, ~~\forall v \in {V}({\cal T}_h).
\en
Namely, $\xi_r$ satisfies
\be \label{resivari}
a(\xi_r, v) = L(v) - a(u_{r-1}, v),   ~~\forall v \in {V}({\cal T}_h),
\en
 which indicates that $\xi_r$ is the residual of the current iterative solution $u_{r-1}$:
\be \label{residefi}
\xi_r = u - u_{r-1}.
\en

 Unfortunately, we can not compute the minimizer $\xi_r$ explicitly, since it would be equivalent to solving the original minimization problem (\ref{minmum}) exactly. Instead, we compute an approximate minimizer $\xi_r^{\theta} \approx \xi_r$ in an appropriate way.

For a positive integer $n_r$ defined later,  when aiming to minimize $J(u_{r-1} + \cdot)$ in the subset $V_{n_r}({\mathcal T}_h) \subset {V}({\cal T}_h)$, a significant complication is that it  is not a linear subspace and it is not closed (topologically) in ${V}({\cal T}_h)$. Hence, even though $J(u_{r-1} + \cdot)$ have an infimum on $V_{n_r}({\mathcal T}_h)$ (since $J(u_{r-1} + )$ has a lower bound zero), there may not be a minimizer in $V_{n_r}({\mathcal T}_h)$. Therefore, one should not aim to completely minimize $J(u_{r-1} + \cdot)$, but instead use a relaxed notion of quasi-minimization as used by \cite{Karniadakis} (for which the existence of an infimum implies the existence of a quasi-minimizer):

 \begin{definition}
 Let $\delta_r$ be a sufficiently small positive parameter, and let $V_{n_r}({\mathcal T}_h)$ be a subset of ${V}({\cal T}_h)$. A function $\xi_r^{\theta} \in V_{n_r}({\mathcal T}_h)$ is said to be a quasi-minimizer of
 $J(u_{r-1} + \cdot)$ associated with $\delta_r$ if the following inequality holds:
 \be \label{iteradisminmum}
 J(u_{r-1} + \xi_r^{\theta}) \leq \mathop{\text{inf}}\limits_{v_r^{\theta}\in V_{n_r}({\mathcal T}_h)} J(u_{r-1} + v_r^{\theta})+\delta_r.
 \en
 \end{definition}
It is clear that the value of the functional $J$ at the quasi-minimizer $\xi_r^{\theta}$ is a good approximation of the infimum.

After solving $\xi_r^{\theta}$, we use it to further get a better approximation $u_r$ to $u$ by setting
\be \label{updateappro}
u_r = u_{r-1} + \xi_r^{\theta}.
\en

The following result demonstrates that $\xi_r^{\theta}\in V_{n_r}({\mathcal T}_h)$ is a good approximation to $\xi_r$ under some mild conditions.

\begin{lemma} \label{residualappro}
Assume that $|||u-u_{r-1}|||>0$, and let $0< \tau_r< 1$ be given. Choose $\delta_r= \frac{\tau^2_r}{2}~ |||u-u_{r-1}|||^2$.
Then there exists $n(\tau_r, u_{r-1})\in \mathbb{N}$ such that, when $n_r\geq n(\tau_r, u_{r-1})$, the quasi-minimizer $\xi_r^{\theta}$ defined in Definition 4.1 has the estimate:
\be \label{redualappro}
||| \xi_r^{\theta} - \xi_r ||| < \tau_r ~|||u-u_{r-1}|||.
\en
\end{lemma}
{\it Proof.}  Noting that the underlying set $V_{n_r}({\mathcal T}_h)$ is not a linear space, we can not use the standard technique to analyze the above estimate, instead, we will use special properties of
quadratic functional $J$.

Define $j(v)= J(u_{r-1} + v)$. Then we have $j'(\xi_r)=0$ and
$$ \langle j''(v)w_1,~w_2\rangle= 2\text{Re}\big\{ a(w_1,~w_2) \big\},\quad \forall v,w_1,w_2 \in V({\cal T}_h). $$
Then, by the second-order Taylor formula, there exists $\xi^{\ast}_r\in
V({\mathcal T}_h)$ such that %(see \cite)
\be \label{tay1}
j(\xi_r^{\theta}) - j(\xi_r)  = \text{Re}\big\{ \langle j''(\xi^{\ast}_r)(\xi_r^{\theta}-\xi_r),\xi_r^{\theta}-\xi_r \rangle \big\} = a(\xi_r^{\theta}-\xi_r,\xi_r^{\theta}-\xi_r)  = ||| \xi_r^{\theta} - \xi_r|||^2,
\en
and
\be\label{tay2}
j(v_r^{\theta}) - j(\xi_r) =  |||v_r^{\theta}-\xi_r |||^2, \quad~\forall v_r^{\theta}\in V_{n_r}({\mathcal T}_h).
\en
Using the quasi-minimization property of $j(\xi_r^{\theta})$, we have
$$ j(\xi_r^{\theta}) - j(\xi_r) \leq \mathop{\text{inf}}\limits_{v_r^{\theta}\in V_{n_r}({\mathcal T}_h)} j(v_r^{\theta}) + \delta_r - j(\xi_r) \leq j(v_r^{\theta}) - j(\xi_r) + \delta_r.$$
Thus, combining (\ref{tay1}) with (\ref{tay2}), we get the quasi-optimal estimate
\be \label{bestappro}
|||\xi_r - \xi_r^{\theta} |||^2 \leq \mathop{\text{inf}}\limits_{v_r^{\theta}\in V_{n_r}({\mathcal T}_h)}|||\xi_r - v_r^{\theta} |||^2 + \delta_r.
\en
Since $|||u-u_{r-1}|||>0$, by the universal approximation property, there exists $n(\tau_r, u_{r-1})\in \mathbb{N}$ and $\tilde{\xi}_r^{\theta}\in  V_{n(\tau_r, u_{r-1})}({\mathcal T}_h)$ such that
\be \label{interappro}
|||\xi_r-\tilde{\xi_r^{\theta}}||| < \frac{\tau_r}{\sqrt{2}}~ |||u-u_{r-1}|||.
\en
Let $n_r\geq n(\tau_r, u_{r-1})$ and choose $v_r^{\theta}=\tilde{\xi_r^{\theta}}$ in (\ref{bestappro}). Substituting (\ref{interappro}) into (\ref{bestappro}) and
noting that $\delta_r= \frac{\tau^2_r}{2}~ |||u-u_{r-1}|||^2$, we get the desired result (\ref{redualappro}).
$\hfill\Box$

The next result addresses the convergence of the proposed recursive algorithm.

\begin{theorem} \label{helmtheom}
Assume that, for $j=1,\cdots,r$, $|||u-u_{j-1}|||>0$,~$0<\tau_j < 1$ and $\delta_j=\frac{\tau^2_j}{2}~ |||u-u_{j-1}|||^2$. Let $n(\tau_j, u_{j-1})\in \mathbb{N}$ and
$\xi_j^{\theta}\in V_{n_j}({\mathcal T}_h)$ with $n_j\geq n(\tau_j, u_{j-1})$ be determined as in Lemma \ref{residualappro} ($j=1,\cdots,r$).
Then the resulting approximate solution $u_r$ has the estimate:
\be
|||u-u_r||| < |||u-u_0|||~\prod_{j=1}^r\tau_j.
\en
\end{theorem}
{\it Proof.} By the iterative scheme (\ref{updateappro}) and the definition (\ref{residefi}) of the residual $\xi_r$, we get
\be \label{initierror}
 |||u-u_r||| = ||| u - ( u_{r-1} + \xi_r^{\theta}) ||| = ||| \xi_r -  \xi_r^{\theta} |||.
 \en
Using the approximation (\ref{redualappro}) recursively, yields
%can directly the result.
\beq \label{iteresti}
|||u-u_r||| & < & \tau_r |||u-u_{r-1}||| \cr
&& \vdots
 \cr
 & < &  |||u-u_0|||~\prod_{j=1}^r \tau_j.
\eq
$\hfill\Box$

%\begin{corollary}
%By the universal approximation property, setting $\tau_r \rightarrow 0$ when $ r \rightarrow \infty$, then we have the strong convergence of a sequence of approximations $u_r$ to the solution to the governing equation as follows.
%\be
%|||u-u_r||| \rightarrow 0 ~~\text{when} ~~r\rightarrow \infty.
%\en
%\end{corollary}

\begin{remark} The key result Lemma 4.2 looks like Proposition 2.5 in \cite{AD}, but Lemma 4.2 was established for a different method from the one introduced in \cite{AD}.
In the proof of Proposition 2.5 in [1], the equality (2.12) in [1] was used (see line-7 in the proof), but this equality does not hold
if the boundedness assumption of the neural network parameters  $\{C_i\}$ is not met (see Remark 2.4 in [1]). This means that the boundedness assumption of the neural network parameters
was implicitly required in Proposition 2.5 of \cite{AD}. As pointed out in Remark 3.1, there is no boundedness requirement on the neural network parameters in the proposed DPWNN method,
so the approximation property Lemma 4.2 cannot be directly obtained by Proposition 2.5 in \cite{AD} and was proved in a different manner from that provided in [1].

%, where the sequence numbers $\{C_i\}$ was required to be bounded.
%the main theoretical results  depend on the assumption of boundedness of sequence numbers $\{C_i\}$. So far, there is no relevant literature that rigorously proves the boundedness of sequence numbers $\{C_i\}$ introduced in %\cite{AD}, even if Figure 6 in \cite{AD} shows that the network parameters are well-bounded and do not exhibit explosive behavior.
\end{remark}

%It has been numerically investigated in \cite[The second-to-last paragraph on p. 17]{ref21} and \cite[Remark 12]{ref11} that, linear parameters $c$ in neural networks may be unbounded, when $n_r$ increases to infinity and $h$ tends to zero.

\begin{remark} The method described in this section is different from that developed in \cite{AD}:
the former is naturally derived from the quasi-minimization problem (\ref{iteradisminmum}), but the latter was derived from a variational problem
in a slightly complex manner. The cost for solving the quasi-minimization problem (\ref{iteradisminmum}) is almost the same as that for solving the complementarity problem (2.16) in \cite{AD},
but there is an extra step to solve the variational problem (2.17) in \cite{AD}. This means that the current method is cheaper than that introduced in \cite{AD}.
%By (\ref{initierror}), we find that, the error of approximate solution $u_r$ is in fact equal to the error of approximate minimizers $\xi_r^{\theta}$ to the residual $\xi_r$ of the iterative solution $u_{r-1}$; while by %Proposition 2.6 of the continuous Galerkin neural network method \cite{AD},  the approximation error of $\varphi_i^{NN}$ to the dual representation $\varphi_i$ directly acts as the convergence factor of the iteration.
\end{remark}

In the following we describe practical algorithms for computing a good approximate solution $u_r$ of $u$.\\
\vspace{0.2cm}
{\bf Algorithm 4.1} (iteratively compute an approximate solution $u_r$)

Give a termination parameter $tol$ and the initial guess $u_0=0$. For $r\geq 1$, let $u_{r-1}\in V({\mathcal T}_h)$ be gotten, and assume that
$J(u_{r-1}) \geq tol$. The approximation $u_r$ is generated as follows:

{\bf Step 1} Choose a proper positive integer $n_r$ and solve the quasi-minimization problem (\ref{iteradisminmum}) by {\bf Algorithms 4.2-4.3} to get $\xi_r^{\theta}\in V_{n_r}({\mathcal T}_h)$;

{\bf Step 2} Define $u_r=u_{r-1}+\xi_r^{\theta}$.

If $J(u_{r})<tol$, then the iteration is terminated. \\

The quasi-minimization problem (\ref{iteradisminmum}) can not be solved directly. Instead, we have to alternatively compute plane wave direction angles and activation coefficients by iterative
algorithms.\\
\vspace{0.2cm}
{\bf Algorithm 4.2} (alternatively compute activation coefficients and direction angles)

Consider the quasi-minimization problem (\ref{iteradisminmum}) with known $u_{r-1}$ and $n_r$. Give a small parameter $\rho>0$ and initial direction angles $\Phi^0$. For known direction angles $\Phi^l$ with $l\geq 0$,
let $V_{\Phi^l}({\mathcal T}_h)$ denote the space $V_{n_r}({\mathcal T}_h)$ with the given direction angles $\Phi^l$.

{\bf Step 1} Solve the following minimization problem by the standard Galerkin method
\be \label{dglsq}
 J(u_{r-1} + \xi^{l}_r)=\min\limits_{{v}\in V_{\Phi^{l}}({\mathcal T}_h)}~J(u_{r-1} + v)
 \en
to get the activation coefficients $c^l$, which are the coordinate vector of the function ${\xi^l_r}$ under the plane wave basis functions with the direction angles $\Phi^l$;

{\bf Step 2} Solve the following quasi-minimization problem by a gradient descent algorithm (see Algorithm 4.3)
\be \label{dggs}
 J(u_{r-1} + \eta^{l+1}_r) \approx \inf\limits_{{v}\in V_{c^l}({\mathcal T}_h)}~J(u_{r-1} + v)
 \en
to get updated direction angles $\Phi^{l+1}$ (which are the direction angles of the plane wave basis functions defining the function $\eta^{l+1}_r$)
%\textcolor{red}{from the old direction angles $\Phi^{l}$,}
where $V_{c^l}({\mathcal T}_h)$ denotes the set $V_{n_r}({\mathcal T}_h)$ with the given activation coefficients $c^l$.

If $||\Phi^{l+1} - \Phi^{l}||_{\infty} < \rho$, then the iteration is terminated; %\textcolor{red}{}

else $l=l+1$ and go to {\bf Step 1}. \\

 A pseudo {\bf Algorithm 1} combining {\bf Algorithm 4.1} and {\bf Algorithm 4.2} summarizes this recursive approach.

\begin{algorithm}[H]
\caption{Discontinuous Plane Wave Neural Network}

\quad Set $u_0=0$ and tolerance $tol > 0$.

\quad {\bf for} $r=1 : \text{maxit}$ {\bf do}

%\quad \quad  $(\xi_{i}^{NN}, J(u_{r-1})) \leftarrow $ PWRes($u_{r-1}, a, L$)

\quad\quad  Given $u_{r-1}$, $n_r$ and a small parameter $\rho>0$.

\quad\quad  Initialize hidden parameters $\Phi^0\in \mathbb{R}^{n_r}$.

\quad\quad  {\bf For} $l = 0 : \text{traincount}$

\quad \quad \quad Compute $c^l$ by {\bf Function} DLSQ-R$(\Phi^l, \sigma, a, L(\cdot)-a(u_{r-1}, \cdot))$.

\quad \quad  \quad Update the propagation wave directions for $\Phi^l$ by Algorithm 4.3.

\quad \quad  \quad {\bf if } $||\Phi^{l+1} - \Phi^{l}||_{\infty} < \rho$

\quad \quad  \quad \quad return $\xi_r^l$ and $J(u_{r-1}+\xi_r^l)$

\quad \quad  \quad {\bf end if}

\quad\quad  {\bf end for}

\quad \quad {\bf if } $(J(u_{r-1}+\xi_r^l) < tol)$

\quad \quad \quad break

\quad \quad {\bf else}

\quad \quad \quad $u_{r} = u_{r-1} +\xi_r^l$

\quad \quad \quad $r = r+1$

\quad \quad {\bf end if}

\quad  {\bf end for}

\quad Return $N=r-1$ and $u_N$.
\end{algorithm}

Here {\bf Function} DLSQ-R$(\Phi^l, \sigma, a, L(\cdot)-a(u_{r-1}, \cdot))$ computes the expansion coefficients of the projection of the error $u-u_{r-1}$ onto the space $V_{n_r}^{h}$ with respect to $a$ and the residue functional $\langle \xi(u_{r-1}), v\rangle = L(v) - a(u_{r-1}, v)$ provided that hidden parameters $\Phi^l$ are given; please see the subsequent section 4.3.

\iffalse
\quad\quad  Initialize hidden parameters $\Phi^{(r)}\in \mathbb{R}^{n_r}$.

\quad \quad Compute corresponding activation coefficients:

\quad \quad $c^{(r)}$ by {\bf Function} dGLSQ-R$(\Phi^{(r)}, \sigma, a, L(\cdot)-a(u_{r-1}, \cdot))$.

\quad \quad {\bf for} each training epoch {\bf do}

\quad \quad  \quad Compute the update propagation wave directions for $\Phi^{(r)}$ by {\bf Algorithm 4.3}.

%Remove
\quad \quad  \quad $\{\zeta_m\}$ corrections for 3D case: Correct associated angles that satisfy $sin(\zeta_m)\approx 0$ by

\quad \quad  \quad adding suitable disturbance term $\varepsilon_{disturb}$ to guarantee the linear independence of the

\quad \quad  \quad  resulting basis functions $\{e^{\text{i}\omega {\bf d}_{l}\cdot{\bf x}}\}$ for different $\vartheta_t$.%^{(k)}

\quad \quad  \quad Compute corresponding activation coefficients:

\quad \quad  \quad $c^{(r)}$ by {\bf Function} dGLSQ-R$(\Phi^{(r)}, \sigma, a, L(\cdot)-a(u_{r-1}, \cdot))$.

\quad \quad  {\bf end for}

\quad \quad  Return $\xi_r^{\theta}(\theta^{(r)})$ and $J(u_{r-1})$.

\fi

\subsection{Discussions on Algorithm 4.1}

%\textcolor{red}{how to choose $h$ and $C_i$ practically}

We address that each approximate minimizer $\xi_r^{\theta} \approx \xi_r$ may be realized in the discontinuous plane wave neural network $V_{n_r}({\mathcal T}_h)$ with a single hidden layer by learning the quasi-minimizer of the loss $J(u_{r-1} + v^{\theta})$ on the set $V_{n_r}({\mathcal T}_h)$.

%Owing to the employed strategy of $hp-$refinement, to achieve higher approximation accuracy, the requirement above on the width $n$ of the network for the proposed discontinuous Galerkin neural network needs to be appropriately increased.
If the width $n_r$ is fixed as $r$ increases, the rate of convergence will become very slow due to the decreased ability of a fixed-width network $V_{n}({\mathcal T}_h)$ to capture higher resolution features of the error as $r$ increases (see \cite[section 2.3.2]{AD}). Thus, we necessarily require that $n_r=n(\tau, u_{r-1})>n_{r-1}=n(\tau, u_{r-2})$ such that
$\xi_r^{\theta}\in V_{n_r}({\mathcal T}_h)$ can capture the higher resolution features of $\xi_r$.

On the other hand, by the $hp$-convergence (\ref{approx}) of the aforementioned discontinuous Galerkin method, the errors of the approximations $u_r$ with refining mehswidth $h$ will
decrease more fast than that with a fixed $h$ when the wave number $\omega$ increases.
 That is, a judicious $hp-$refinement strategy will be the most attractive option to solve the wave problem with large wave numbers.

Numerical results in section 6 will validate that, when choosing $h\approx \mathcal{O}(\frac{4\pi}{\omega})$ and gradually increasing the width $n_r$ of the network at each iteration in Algorithm 4.1,
the resulting approximations $u_r$ can reach a given accuracy. In fact, the iteration counts of Algorithm 4.1 are independent of wave numbers when $n_r$ is chosen in a proper manner.
In particular, only no more than ten iterations of Algorithm 4.1 may guarantee the convergence for the tested examples.

\subsection{Calculation of the coefficients $c^l$}
Note that, when the direction angle parameters ${\Phi^{l}}$ are known, the discontinuous plane wave neural network set $V_{n_r}({\mathcal T}_h)$ constitutes a linear space, consisting of $N\times n_r$ neural network basis functions $\{\psi_{kj}\}$, which satisfy
\begin{eqnarray}
\psi_{kj}(x)=\left\{\begin{array}{ll}
 e^{\text{i} \omega {\bf d}_j\cdot {\bf x}},~~{\bf x}\in\Omega_k,\\
 0,~~{\bf x}\in\Omega_l~~\mbox{satisfying}~~l\neq k
\end{array}\right.~~(k,l=1,\cdots,N;~j=1,\cdots, n_r).
\label{neubasis}
\end{eqnarray}
 Then a PWLS approximation associated with (\ref{dglsq}) can be computed from $V_{n_r}({\mathcal T}_h)$ by choosing the parameters $c^{l}$ as follows:

\begin{eqnarray}
\left\{\begin{array}{ll}  {\cal A}c^{l} = b, \\
{\cal A}^{k,l}_{s,j} = a(\psi_{lj},\psi_{ks}), ~~b^k_s=L(\psi_{ks}) - a(u_{r-1},\psi_{ks}).
\end{array}\right.
\label{galerlsq}
\end{eqnarray}
This procedure defines a function which computes the expansion coefficients of the projection of $u-u_{r-1}$ onto the space $V_{n_r}({\mathcal T}_h)$ with respect to $a$ and the residual $\langle \xi(u_{r-1}), v\rangle=L(v) - a(u_{r-1}, v)$.

Once each training epoch in Algorithm 4.2 is terminated, we can assemble the ultimate  basis functions $\{\psi_{kj}\}$ associated with the direction angle parameters ${\Phi^{l}}$ and the corresponding coefficient $c^{l}$ in a proper way, and finally obtain the current approximate minimizer $\xi_r^{\theta}$ of the loss $J(u_{r-1} + v_r^{\theta})$ on the set $V_{n_r}({\mathcal T}_h)$.

\subsection{Calculation of the direction angles} \label{hlempara}
To calculate optima of the objective functional (\ref{iteradisminmum}) associated with discontinuous plane wave neural networks, the modified version of a gradient-based Adam optimizer \cite{luo} with full batch is employed to update parameter set ${\Phi}$ of hidden direction angles in (\ref{dggs}) at each training epoch. In the first training step, the direction angle set ${\Phi}^0$ is initialized uniformly, which will be described in the numerical test. In later iteration steps, the direction angle set is updated via an randomly shuffled Adam optimizer \cite{luo}. In particular, a single step of the Adam optimizer is adopted and the set of $N_G$ quadrature nodes is selected as training data to calculate the loss function.

The calculation process of updated direction angle set $\Phi^{l+1}$ can be described in Algorithm 4.3.  \\
\vspace{0.2cm}
%\begin{algorithm}[H]
{\bf Algorithm 4.3} (compute direction angles with given activation coefficients)

\quad Initialize $\Phi_{1,0}=\Phi^0$, $m_{1,-1}=\nabla J(\Phi^0)$ and $v_{1,-1} = \text{max}_i \nabla f_i(\Phi^0) \circ \nabla f_i(\Phi^0)$.

\quad {\bf for} $s=1 : \infty$ {\bf do}

\quad \quad Sample $\{\tau_{s,0},\tau_{s,1},\cdots, \tau_{s,N_G-1} \}$ as a random permutation of $\{0,1,2,\cdots, N_G-1\}$.

\quad \quad  {\bf for} $i=0 : N_G-1$ {\bf do}

\quad \quad  \quad $m_{s,i}=\beta_1 m_{s,i-1} + (1-\beta_1)\nabla f_{\tau_{s,i}}(\Phi_{s,i})$.

\quad \quad  \quad $v_{s,i}=\beta_2 v_{s,i-1} + (1-\beta_2) \nabla f_{\tau_{s,i}}(\Phi_{s,i}) \circ  \nabla f_{\tau_{s,i}}(\Phi_{s,i})$.

\quad \quad  \quad $\Phi_{s,i+1}= \Phi_{s,i} - \frac{\eta_s}{\sqrt{v_{s,i}}+\epsilon} \circ m_{s,i}$.

\quad \quad  {\bf end for}

\quad \quad  $\Phi_{s+1,0}=\Phi_{s,N_G}, ~v_{s+1,-1}=v_{s,N_G-1}, ~m_{s+1,-1}=m_{s,N_G-1}$.

\quad  {\bf end for}
%\end{algorithm}

%\begin{remark}\label{diffongrad}
% We would like to emphasize the reason why the above Theorem holds is that the sesquilinear form $a$ induced by the PWLS method is bounded, Hermitian positive definite, and the semilinear form $L$ is bounded, which complies with the de facto conditions of \cite[Section 2.1]{AD}. Of course, the introduction of the intrinsic PWLS method brings significantly different definitions of bilinear forms $a$, which further induces the completely different definitions of norms $|||\cdot|||$. Additionally, due to the injection of complex Euclidean space in the considered model and plane wave activation functions, the residual functional $r(u_{r-1})$ defined by (\ref{resiformulaiter}) is the module of the pre-residue $\langle \xi(u_{r-1}), v\rangle =L(v) - a(u_{r-1}, v)$, which is a visible another difference from that of \cite{AD}.
%It leads to the further difference on the computation of derivative of $\eta(u_{r-1},v)$ with respect to the training parameters $D$.
% \end{remark}

In Algorithm 4.3, $\circ$ denotes the component-wise product, the division and square-root operator are component-wise as well. The parameter $m$ denotes the 1st-order momentum and $v$ denotes the 2nd-order momentum. They are weighted averaged by hyperparameter $\beta_1, \beta_2$, respectively. Denote $\Phi_{s,i}, m_{s,i}, v_{s,i}$ as the value of $\Phi, m, v$ at the $s-$th outer loop (epoch) and $i-$th inner loop (batch), respectively. As suggested in the article \cite{luo}, we set $\beta_1=0.9, ~\beta_2=0.999$, choose $\eta_s=\frac{\eta_1}{\sqrt{kN_G}}$ as a learning rate scheduler to reduce the learning rate during the process, and set $\epsilon=10^{-8}$.

%Particularly, a single step of the Adam optimizer is adopted and the set of $N_G$ quadrature nodes is selected as training data to calculate the loss function.

%As aforementioned in Remark \ref{diffongrad},

We would like to describe some details for computation of gradient of $J(u_{r-1} + v^{\theta})$ with respect to the training variables $\Phi$, which will be used to update the direction angles. % in Algorithm 4.3 stated above.

%an randomly shuffled Adam optimizer \cite{luo}.

 %involving the calculation of derivative $\nabla_D \eta(u_{r-1},v)$ of the quantity $\eta(u_{r-1},v)$ defined in (\ref{etaquan}).

By the direct calculation, we have
\be
\nabla_{\Phi} J(u_{r-1} + v^{\theta}) = \sum_{k=1}^N\int _{\gamma_k} \mathbb{F}(\mathbb{J}_b) ds
 + \sum_{j\not=k} \bigg( \alpha \int_{\Gamma_{kj}} \mathbb{F}(\mathbb{J}_{\alpha})
 ds + \beta \int_{\Gamma_{kj}} \mathbb{F}(\mathbb{J}_{\beta})  ds \bigg), \label{4.new}
\en
where $\mathbb{F}(\mathbb{J}) = 2\text{Re}\big( \overline{\mathbb{J}} ~ \nabla_{\Phi}\mathbb{J} \big)$ for any scalar quantity $\mathbb{J}$;
$$ \mathbb{J}_b= (\partial_{\bf n}+i\omega )(u_{r-1} + v^{\theta})-g,\quad \mathbb{J}_{\alpha} = (u_{r-1} + v^{\theta})|_{\Omega_k}-(u_{r-1} + v^{\theta})|_{\Omega_j}$$
and
$$\mathbb{J}_{\beta} = \partial_{\text{\bf n}_k}(u_{r-1} + v^{\theta})+\partial_{\text{\bf n}_j}(u_{r-1} + v^{\theta}). $$

%\be
%\nabla_D \eta(u_{r-1},v) = \frac{1}{\eta(u_{r-1},v)} \bigg(\text{Re}\big(\xi(u_{r-1}, v)\big) \cdot \text{Re}\big(\xi(u_{r-1}, \nabla_D v)\big) +  \text{Im}\big(\xi(u_{r-1}, v)\big) \cdot \text{Im}\big(\xi(u_{r-1}, \nabla_D v)\big) \bigg).
%\en

The formula (\ref{4.new}) can be written more clearly for a special form of $v^{\theta}$. In two dimensional case, taking
$v^{\theta} =  \sum_{l=1}^n c_l^{(k)} \text{e}^{\text{i} \omega {\bf d}_l^{(k)} \cdot {\bf x}}, ~{\bf x}=(x_1, x_2)^T \in \Omega_k \in {\cal T}_h$, we get $\nabla_{\Phi^{(k)}} v^{\theta} =(\partial_{\vartheta_1^{(k)}}v^{\theta}, \cdots, \partial_{\vartheta_n^{(k)}}v^{\theta})^T$, where
\be \nonumber
\partial_{\vartheta_m^{(k)}}v^{\theta} = i\omega  c_m^{(k)} \big(-x_1 sin(\vartheta_m^{(k)}) + x_2 cos(\vartheta_m^{(k)})\big)~e^{\text{i}\omega {\bf d}_{m}^{(k)}\cdot {\bf x}}.
\en
In three dimensional case, taking $v^{\theta} =  \sum_{m=1}^{m^{\ast}}\sum_{t=1}^{t^{\ast}}  c_l^{(k)} ~e^{\text{i}\omega
{\bf d}_{l}^{(k)}\cdot {\bf x}}, ~l=(m-1)k^{\ast}+t, ~{\bf x}=(x_1, x_2, x_3)^T \in \Omega_k \in {\cal T}_h$, we get $\nabla_{\Phi^{(k)}} v^{\theta} = (\partial_{\vartheta_1^{(k)}}v^{\theta}, \cdots, \partial_{\vartheta^{(k)}_{t^{\ast}}}v^{\theta}, \partial_{\zeta_1^{(k)}}v^{\theta}, \cdots, \partial_{\zeta^{(k)}_{m^{\ast}}}v^{\theta})^T$, where
\be \nonumber
\partial_{\vartheta_t^{(k)}}v^{\theta} = i\omega \sum_{m=1}^{m^{\ast}} c_l^{(k)} \big(-x_1 sin(\zeta_m^{(k)})sin(\vartheta_t^{(k)}) + x_2 sin(\zeta_m^{(k)})cos(\vartheta_t^{(k)})\big)~e^{\text{i}\omega
{\bf d}_{l}^{(k)}\cdot {\bf {\bf x}}},
\en
and
\be\nonumber
\partial_{\zeta_m^{(k)}}v^{\theta} = i\omega \sum_{t=1}^{t^{\ast}} c_l^{(k)}  \big(x_1 cos(\zeta_m^{(k)})cos(\vartheta_t^{(k)}) + x_2 cos(\zeta_m^{(k)})sin(\vartheta_t^{(k)})-x_3sin(\zeta_m^{(k)})\big)~e^{\text{i}\omega
{\bf d}_{l}^{(k)}\cdot {\bf x}}.
\en

\section{A DPWNN method for Maxwell's equations}
In this section, we extend the DPWNN method developed in the last section to the second order homogeneous Maxwell equations with the lowest order absorbing boundary condition in three space dimensions:
\begin{equation} \label{maxeq}
\left\{ \begin{aligned}
     &  \nabla\times(\frac{1}{i\omega\mu}\nabla\times {\bf E}) + i\omega\varepsilon {\bf E} = 0 & \text{in}\quad
\Omega,\\
    & -{\bf E}\times {\bf n}+\frac{\sigma}{i\omega\mu}((\nabla\times {\bf E})\times  {\bf n})\times {\bf n}={\bf g} & \text{on}\quad\gamma.
                          \end{aligned} \right.
                          \end{equation}
Here $\omega>0$ is the temporal frequency of the field, and $
{\bf g}\in {\bf L}_T^2(\partial\Omega)^3$. The material coefficients
$\varepsilon,\mu$ and $\sigma$ are assumed to be piecewise constant in the
whole domain. In particular, if $\varepsilon$ takes complex valued,
then the material is an absorbing medium; else the material is a non-absorbing medium (see \cite{hmm}).

\subsection{A minimization problem}
For an element $\Omega_k$, let $\mbox{
\bf H}(\text{curl};\Omega_k)$ denote the standard Sobolev space. Set
\begin{equation}
 {\bf V}(\Omega_k)=\bigg\{{\bf E}_k\in \mbox{
\bf H}(\text{curl};\Omega_k); ~ {\bf E}_k \text{~elementally satisfies the
first equation of (\ref{maxeq})}\bigg\}.
\end{equation}
Define
 \be
{\bf V}({\cal T}_h) = \{ {\bf v}\in {\bf L}^2(\Omega): ~{\bf v}|_{\Omega_k}\in {\bf V}(\Omega_k) ~\forall \Omega_k\in {\cal T}_h  \}.
\en

 For each local interface $\Gamma_{kj}$ $(k<j)$, we define the jump on
$\Gamma_{kj}$ as follows (note that ${\bf n}_k=-{\bf n}_j$):
\be
\llbracket {\bf F}\times  {\bf n} \rrbracket  =  {\bf F}_k\times {\bf n}_k + {\bf F}_j\times {\bf n}_j.
 \en

Then the reference problem (\ref{maxeq}) will be solved by finding the local electric field ${\bf E}_k$ such that
\begin{equation} \label{elementmaxeq}
\left\{ \begin{aligned}
     &  \nabla\times(\frac{1}{i\omega\mu}\nabla\times {\bf E}_k) + i\omega\varepsilon {\bf E}_k = 0 & \text{in}\quad
\Omega_k,\\
    & -{\bf E}_k\times {\bf n}+\frac{\sigma}{i\omega\mu}((\nabla\times {\bf E}_k)\times  {\bf n})\times {\bf n}={\bf g} & \text{on}\quad\gamma_k,
                          \end{aligned} \right.
                          \end{equation}
with the transmission conditions on each interface $\Gamma_{kj}$
\be  \label{tangjump}
\llbracket{ {\bf E}}\times{ {\bf n}} \rrbracket = 0\quad \text{and} \quad \llbracket {1\over i\omega\mu}(\nabla\times{{\bf E}})\times {\bf n} \rrbracket = 0.
\en

Based on this observation, we define the functional (refer to \cite{hy2})
\beq \label{max_mini_func}
J({\bf F})&=&\sum_{k=1}^N\int _{\gamma_k}| - {\bf F}_k\times  {\bf n}_k
 + \frac{\sigma}{i\omega\mu}((\nabla\times  {\bf F}_k)\times
{\bf n}_k) \times {\bf n}_k - {\bf g}|^2~ds \cr&+&
 \sum_{k<j} \Bigg(\int_{\Gamma_{kj}}
\rho_1 |\llbracket{ {\bf F}}\times{ {\bf n}} \rrbracket|^2~~ds
+
 \int_{\Gamma_{kj}} \rho_2 |\llbracket {1\over i\omega\mu}(\nabla\times{{\bf F}})\times {\bf n} \rrbracket|^2 ~ds \bigg), \quad {\bf F}\in {\bf V}({\cal T}_h),
\eq
where $\rho_1$ and $\rho_2$ are two fixed positive constants.

Consider the minimization problem: find ${\bf E}\in {\bf V}({\cal T}_h)$ such that
\begin{equation} \label{maxminmum}
J({{\bf E}})=\min\limits_{{{\bf F}}\in {\bf V}({\cal T}_h)}~J({\bf F}).
\end{equation}

The PWLS variational problem associated with the minimization problem (\ref{maxminmum}) can be expressed as follows (see \cite{hy2}): find ${\bf E}\in {\bf V}({\cal T}_h)$ such that, for $\forall ~ {\bf F}\in {\bf V}({\cal T}_h)$,
\beq
%\begin{equation}
%\begin{split}
 &~~{\small \sum_{k=1}^N\int _{\gamma_k}\bigg(-{\bf E}_k\times  {\bf n}_k + \frac{\sigma}{i\omega\mu}\big((\nabla\times  {\bf E}_k)\times
{\bf n}_k\big) \times {\bf n}_k - {\bf g}\bigg)\cdot\overline{-{\bf F}_k\times {\bf n}_k + \frac{\sigma}{i\omega\mu}\big((\nabla\times  {\bf F}_k)\times
{\bf n}_k\big) \times {\bf n}_k}~ds}
\cr
&\quad\quad+
 \sum_{k<j}\int_{\Gamma_{kj}}\Bigg(\rho_1 \llbracket  {\bf E}\times  {\bf n}
  \rrbracket \cdot \overline{\llbracket  {\bf F}\times  {\bf n} \rrbracket}
+ \rho_2
\llbracket {1\over i\omega\mu}(\nabla\times {\bf E})\times {\bf n} \rrbracket \cdot \overline{ \llbracket {1\over i\omega\mu}(\nabla\times {\bf F})\times {\bf n} \rrbracket }\bigg) ~ds =0.\label{MAXVAR}
%\end{split}
%\label{MAXVAR}
%\end{equation}
\eq

Define the sesquilinear form $a(\cdot, \cdot)$ by
\beq \nonumber
 a({\bf E},{\bf F}) =  \sum_{k<j}\int_{\Gamma_{kj}}\Bigg(\rho_1 \llbracket  {\bf E}\times  {\bf n}
  \rrbracket \cdot \overline{\llbracket  {\bf F}\times  {\bf n} \rrbracket}
+ \rho_2
\llbracket {1\over i\omega\mu}(\nabla\times {\bf E})\times {\bf n} \rrbracket \cdot \overline{ \llbracket {1\over i\omega\mu}(\nabla\times {\bf F})\times {\bf n} \rrbracket }\bigg) ~ds
 \cr
 +\sum_{k=1}^N\int _{\gamma_k} \bigg(-{\bf E}_k\times  {\bf n}_k + \frac{\sigma}{i\omega\mu}\big((\nabla\times  {\bf E}_k)\times
{\bf n}_k\big) \times {\bf n}_k \bigg) \cdot \overline{-{\bf F}_k\times {\bf n}_k + \frac{\sigma}{i\omega\mu}\big((\nabla\times  {\bf F}_k)\times
{\bf n}_k\big) \times {\bf n}_k}~ds,
 %~\forall v\in V({\mathcal T}_h),
\eq
and semilinear form $L({\bf F})$ by
\be \nonumber
L({\bf F})=\sum_{k=1}^N\int _{\gamma_k} {\bf g} \cdot ~\overline{-{\bf F}_k\times {\bf n}_k + \frac{\sigma}{i\omega\mu}\big((\nabla\times  {\bf F}_k)\times
{\bf n}_k \big) \times {\bf n}_k }~ds.
\en
 Then (\ref{MAXVAR}) can be written as the following variational equation:
 \begin{eqnarray}\label{maxpwlsvar}
\left\{\begin{array}{ll} \text{Find} ~~{\bf E}\in {\bf V}({\cal T}_h) ~~s.t. \\
a({\bf E},{\bf F})=L({\bf F}) \quad \forall {\bf F}\in {\bf V}({\cal T}_h).
\end{array}\right.
\end{eqnarray}

Obviously, $a({\bf F},{\bf F})$ is a norm on $V({\mathcal T}_h)$, which can be denoted by $||| {\bf F} |||^2$. Besides, by the Cauchy-Schwarz inequality, $a(\cdot,\cdot)$ satisfies the boundedness and $L({\bf F})$ is bounded.  Furthermore, a standard application of the Risez representation theorem \cite{FO} implies the existence of a unique ${\bf E}\in {\bf V}({\cal T}_h)$ satisfying (\ref{maxpwlsvar}).

\begin{remark} We would like to give some explanations to the considered PWLS method. From the viewpoint of mathematics, all the DG-type method for Maxwell equations should follow the basic rule: the resulting solution should keep the original tangential continuity of the electric field and its curl across all the element interfaces (see (\ref{tangjump})). In theory there are different ways to meet this rule. The PWLS method may be the simplest DG-type method, which was completely defined by the jumps of tangential traces of the electric field and its curl (see (\ref{MAXVAR})), but the other DG-type methods involve various numerical fluxes (for example, Lax-Friedrich flux), comparing \cite{hy2} with \cite{pwdg}.

%We would like to point out that the choice of the jump quantity defined in (\ref{tangjump}) is to ensure the tangential continuity of the electric and magnetic fields required by the Equation (\ref{maxeq}). It has no relation to %the numerical flux defined by the discontinuous Galerkin method; see the PWDG method \cite[section 3]{pwdg} and the least squares methods \cite {mhp}. Since our method originates from a defined quadratic functional, thus %provided two parameters $\alpha$ and $\beta$ are chosen as positive numbers, the variational formulation (\ref{maxpwlsvar}) is theoretically stable regardless of the wavenumber values. Besides, as pointed out in the theoretical %analysis of \cite[Theorem 4.3]{hy2} and \cite[Theorem 5.7]{pwdg}, the $p-$version of the plane wave methods, also called the spectral approach, is immune to the pollution effect for the case of large wavenumbers. Related %numerical validations can also refer to \cite[section 4]{ref21}.
\end{remark}

\subsection{Discontinuous plane wave neural networks}

 At first, a suitable family of plane waves, solutions of the constant coefficient Maxwell equations, are generated  on every $\Omega_k \in {\cal T}_h$ by choosing $n$ unit propagation directions ${\bf d}_l^{(k)} ~(l=1,\cdots,n)$ (similarly to the Helmholtz case, we use the Algorithm 4.2 to produce them), and then defining every unit real polarization vector ${\bf g}_l^{(k)}$ orthogonal to ${\bf d}_l^{(k)}$. The associated
complex polarization vectors ${\bf f}_l^{(k)}$ and ${\bf f}_{l+n}^{(k)}$ are defined by, via the propagation directions and polarization vectors,
\be \label{polarizavec}
{\bf f}_{l}^{(k)}={\bf g}_{l}^{(k)}, \quad {\bf f}_{l+n}^{(k)}={\bf g}_{l}^{(k)}\times
{\bf d}_{l}^{(k)}~~~(l=1,\cdots,n).
\en

We define two sets of functions ${\bf E}_{l}^{(k)}$ (refer to \cite{pwdg}):
\begin{equation}
{\bf E}_{l}^{(k)}=\sqrt{\mu}~{\bf f}_{l}^{(k)}~e^{\text{i}\kappa
{\bf d}_{l}^{(k)}\cdot{\bf x}}, \quad
{\bf E}_{l+n}^{(k)}=\sqrt{\mu}~ {\bf f}_{l+n}^{(k)}~e^{\text{i}\kappa
{\bf d}_{l}^{(k)}\cdot{\bf x}}~~~
(l=1,\cdots, n),\label{maxeq17}
\end{equation}
where $\kappa=\omega\sqrt{\mu\varepsilon}$. It is easy to verify that two pairs of  functions ${\bf E}_{l}^{(k)}$ and ${\bf E}_{l+n}^{(k)}$ satisfy the elemental Maxwell system (\ref{maxeq}), respectively.

Similarly to subsection \ref{helmcase}, plane wave propagation directions  $\{{\bf d}_{l}^{(k)}\}$ based on spherical coordinates are defined by (\ref{3ddirec}). For a given ${\bf d}_l^{(k)}$, the vectors ${\bf f}_l^{(k)}$ that are orthogonal to ${\bf d}_l^{(k)}$ are not unique. In the following, when we write ${\bf d}_l^{(k)}=(a_l  ~~b_l  ~~c_l)^T$ (which satisfies $a_l^2 +b_l^2 +c_l^2=1$), and assuming that $|b_l|<1$, we always choose ${\bf g}_l^{(k)}$ as follows
$${\bf g}_l^{(k)} = \bigg(\frac{a_lb_l}{\sqrt{1-b_l^2}} ~~~ -\sqrt{1-b_l^2}  ~~~ \frac{-a_l^2b_l+b_l(1-b_l^2)}{c_l\sqrt{1-b_l^2}} \bigg)^T.$$

Subsequently, a discontinuous vector plane wave neural network consists of a single hidden layer of $2n\in \mathbb{R}$ neurons on each element $\Omega_k\in {\cal T}_h$, defining an output function ${\bf E}^{\theta}|_{\Omega_k}: \mathbb{R}^3 \rightarrow \mathbb{C}$ as follows:
\be \label{max_nn_app}
{\bf E}^{\theta}(x)|_{\Omega_k} = \sum_{m=1}^{m^{\ast}}\sum_{t=1}^{t^{\ast}}
\big( c_l^{(k)} {\bf E}_{l}^{(k)} + c_{l+n}^{(k)} {\bf E}_{l+n}^{(k)} \big), ~l=(m-1)k^{\ast}+t.
\en

For simplicity, we will use the elementwise notation $\zeta|_{\Omega_k}=[\zeta_1^{(k)}, \cdots, \zeta^{(k)}_{m^{\ast}}]^T \in \mathbb{R}^{m^{\ast}}$, $\vartheta|_{\Omega_k}=[\vartheta_1^{(k)}, \cdots, \vartheta^{(k)}_{t^{\ast}}]^T \in \mathbb{R}^{t^{\ast}}$, and $c|_{\Omega_k}=[c^{(k)}_1, \cdots, c^{(k)}_{2n}]^T\in \mathbb{C}^{2n}$. A set of these parameters over a collection of all elements is denoted by $\theta=\{\vartheta;\zeta; c\}$. In particular, a set of plane wave propagation angles defined on every element is collectively denoted by $\Phi=\{\vartheta;\zeta\}$.

%A set of nonlinear and linear parameters over a collection of all elements is denoted by $\theta=\{\vartheta;\zeta; c\}$. Particularly, a set of nonlinear parameters defined on every element, is collectively denoted by $\Phi=\{\vartheta;\zeta\}$.

Denote ${\bf V}_n({\cal T}_h)$ as the set of all functions of the form (\ref{max_nn_app}), which is defined by
\be \label{maxvsign}
{\bf V}_n({\cal T}_h) := \bigg\{{\bf F}: ~{\bf F}(x)|_{\Omega_k} =  \sum_{m=1}^{m^{\ast}}\sum_{t=1}^{t^{\ast}}
\big( c_l^{(k)} {\bf E}_{l}^{(k)} + c_{l+n}^{(k)} {\bf E}_{l+n}^{(k)} \big), ~l=(m-1)k^{\ast}+t,~\forall \Omega_k\in {\cal T}_h  \bigg\}.
\en

\iffalse
Given bounded constant $C>0$, we define
\be\label{MAXbounnn}
V_{n,C}^{h}:= \bigg\{  F \in V_n^{h}({\cal T}_h): ||\theta||\leq C   \bigg\},
\en
where $||\theta_k||:=\text{max}|\zeta_m|+\text{max}|\vartheta_k|+\text{max}|c_j|$ for each element $\Omega_k\in {\cal T}_h$,  and $||\theta||=\text{max}_k||\theta_k||$.
\fi

\subsection{A discontinuous plane wave neural network method}
 As in section \ref{helmdis}, we can design iterative algorithms to generate element-wise plane wave neural network solutions approximately satisfying the minimization problem (\ref{maxminmum}). The following result can be proved as in the proof of Theorem \ref{helmtheom}.

 %The next result addresses the convergence of the proposed recursive algorithm.

\begin{theorem} Assume that, for $j=1,\cdots,r$, $|||{\bf E}-{\bf E}_{j-1}|||>0$,~$0<\tau_j < 1$ and $\delta_j=\frac{\tau^2_j}{2}~ |||{\bf E}-{\bf E}_{j-1}|||^2$. Let $n(\tau_j, {\bf E}_{j-1})\in \mathbb{N}$ and
${\bm\xi}_j^{\theta}\in {\bf V}_{n_j}({\mathcal T}_h)$ with $n_j\geq n(\tau_j, {\bf E}_{j-1})$ be determined as in Lemma \ref{residualappro} ($j=1,\cdots,r$).
Then the resulting approximate solution ${\bf E}_r$ has the estimate:
\be
|||{\bf E}-{\bf E}_r||| < |||{\bf E}-{\bf E}_0|||~\prod_{j=1}^r\tau_j.
\en
\end{theorem}

%In this section, we present
%Similarly to subsection \ref{hlempara},
Notice that updating the set $\Phi$ of hidden direction angles requires the calculation of the gradient $\nabla_{\Phi} J({\bf E}_{r-1},{\bf F}^{\theta})$. By the direct reduction, we get

\be
\nabla_{\Phi} J({\bf E}_{r-1} + {\bf F}^{\theta}) = \sum_{k=1}^N\int _{\gamma_k} \mathbb{F}(\mathbb{J}_b) ds
 + \sum_{j\not=k} \bigg( \rho_1 \int_{\Gamma_{kj}} \mathbb{F}(\mathbb{J}_{\rho_1})
 ds + \rho_2 \int_{\Gamma_{kj}} \mathbb{F}(\mathbb{J}_{\rho_2})  ds \bigg),
\en
where $\mathbb{F}(\mathbb{J}) = 2\text{Re}\big(  (\nabla_{\Phi}\mathbb{J})^T ~ \overline{\mathbb{J}} \big)$ for any vector quantity ${\bf \mathbb{ J}}$;
$$ \mathbb{J}_b= \bigg(-({\bf E}_{r-1} + {\bf F}^{\theta})\times  {\bf n} + \frac{\sigma}{i\omega\mu}\big((\nabla\times  ({\bf E}_{r-1} + {\bf F}^{\theta}))\times
{\bf n}\big) \times {\bf n} - {\bf g}\bigg),$$
$$ \mathbb{J}_{\rho_1} = \llbracket  ({\bf E}_{r-1} + {\bf F}^{\theta})\times  {\bf n} \rrbracket,\quad \mathbb{J}_{\rho_2} = \llbracket {1\over i\omega\mu}\big(\nabla\times ({\bf E}_{r-1} + {\bf F}^{\theta}) \big)\times {\bf n} \rrbracket.$$

%\be
%\nabla_D \eta(E_{r-1},F) = \frac{1}{\eta(E_{r-1},F)} \bigg(\text{Re}\big(\xi(E_{r-1}, F)\big) \cdot \text{Re}\big(\xi(E_{r-1}, \nabla_D F)\big) +  \text{Im}\big(\xi(E_{r-1}, F)\big) \cdot \text{Im}\big(\xi(E_{r-1}, \nabla_D F)\big) \bigg).
%\en
Explicitly, taking for example ${\bf F}^{\theta}=  \sum_{m=1}^{m^{\ast}}\sum_{t=1}^{t^{\ast}}
\big( c_l^{(k)} {\bf E}_{l}^{(k)} + c_{l+n}^{(k)} {\bf E}_{l+n}^{(k)} \big), ~l=(m-1)k^{\ast}+t;~~{\bf x}\in \Omega_k \in {\cal T}_h$, we obtain
\beq \nonumber %\label{naDF} %& &
\nabla_{\Phi^{(k)}} {\bf F}^{\theta} = \sum_{m=1}^{m^{\ast}}\sum_{t=1}^{t^{\ast}} \bigg(c_l^{(k)} \sqrt{\mu}~ \big(\nabla_{\Phi^{(k)}} {\bf f}_l^{(k)} + i\kappa {\bf f}_l^{(k)} x^T \nabla_{\Phi^{(k)}} {\bf d}_l^{(k)} \big)
\\
+ c_{l+n}^{(k)} \sqrt{\mu}~  \big(\nabla_{\Phi^{(k)}} {\bf f}_{l+n}^{(k)} + i\kappa {\bf f}_{l+n}^{(k)} x^T \nabla_{\Phi^{(k)}} {\bf d}_l^{(k)} \big) \bigg) ~e^{\text{i}\kappa
{\bf d}_{l}^{(k)}\cdot{\bf x}},
\eq
 where
 $\nabla_{\Phi^{(k)}} {\bf d}_l^{(k)}=(\nabla_{\Phi^{(k)}} a_l ~\nabla_{\Phi^{(k)}} b_l ~\nabla_{\Phi^{(k)}} c_l)^T$, and the calculation of $\nabla_{\Phi^{(k)}} {\bf f}_l^{(k)}$ and $\nabla_{\Phi^{(k)}} {\bf f}_{l+n}^{(k)}$ is similarly to $\nabla_{\Phi^{(k)}} {\bf d}_l^{(k)}$'s. %,  and $\nabla_D W_l \cdot x=x^t \nabla_D W_l$.
 %\end{remark}

%\iffalse

%\input{chapter/chap-Num}
\section{Numerical experiments} \label{numerical}

In this section we report some numerical results to confirm the effectiveness of the proposed DPWNN methods.
%We consider two-dimensional and three-dimensional Helmholtz equation and three-dimensional Maxwell equations.% with large wave numbers.
 For the case of Helmholtz equation in homogeneous media, as pointed out in section 3 of \cite{hy}, we choose the weighted parameters $\alpha$ and $\beta$ in the variational
problem (\ref{pwlsvar}) as $\alpha=\omega^2$ and $\beta=1$.
Besides, as pointed out in \cite{hmm}, %(p.733),
 the permittivity of the homogeneous material in the computational
domain for Maxwell's equations is in general complex valued (i.e.,
the material is an absorbing medium), so we assume that \(\varepsilon=1+i\).
Moreover, we apply the PWLS method with
\(\rho_1=\rho_2=1\) to solve time-harmonic Maxwell's
equations (\ref{maxeq}) in three-dimensional absorbing media.

 For the examples tested in this section, we adopt a uniform triangulation \(\mathcal {T}_h\) for the domain \(\Omega\) as follows. \(\Omega\) is divided into rectangles (for two-dimensional case), or small cubes (for three-dimensional case) with equal meshwidth.
 %, where \(h\) is the length of the longest edge of the elements.

The initial guesses $u_0$ and ${\bf E}_0$ are both uniformly set to be $0$. We %\textcolor{red}{uniformly set the tolerance $tol$ to be $tol = 10^{-6}$,} and
adaptively construct a discontinuous plane wave neural network that achieves an approximation with $J(u_{r-1})< tol$ ($J({\bf E}_{r-1})< tol$ for Maxwell's equations). Moreover, in Algorithm 4.2, the termination condition for training epochs is that $||\nabla_{\Phi} J(u_{r-1} + v^{\theta})||_{\infty} < 10^{-6}$ ($||\nabla_{\Phi} J({\bf E}_{r-1} + {\bf F}^{\theta})||_{\infty} < 10^{-6}$ for Maxwell's equations).%, which will be determined in each numerical test.

%for the Helmholtz equation and $tol = 5.0\times 10^{-5}$ for the Maxwell equation

A set $\Phi$ of plane wave propagation angles introduced in section \ref{hlempara} is uniformly initialized as follows. For the 2D case, the hidden parameters are initialized such that
\be \nonumber
\alpha_j^{(r)}=-\pi+\frac{2\pi}{n_r}j, ~j=1,\cdots,n_r.
\en
 For the 3D case, the hidden parameters are initialized such that
\be  \nonumber
\zeta_m^{(r)}= \frac{\pi}{m^{\ast}_{r}-1}(m-1)+\frac{\pi}{3m^{\ast}_{r}}, ~m=1,\cdots,m^{\ast}_{r}
 \en
 and
 \be  \nonumber
 \vartheta_t^{(r)}=-\pi+\frac{2\pi}{t^{\ast}_{r}}t, ~t=1,\cdots,t^{\ast}_{r}=2m^{\ast}_{r},
 \en
where $m^{\ast}_{r}$ and $t^{\ast}_{r}$ denote the widths $m^{\ast}$ and $t^{\ast}$ (defined in (\ref{3ddirec})) of the network architecture for quasi-minimization iteration $r$, respectively.

To evaluate the objective function and train the discontinuous plane wave neural networks, we employ a fixed complex Gauss-Legendre quadrature rule to approximate all inner products. We perform all computations on a Dell Precision T5500 graphics workstation (2*Intel Xeon X5650 and 6*12GECC) using homemade \textsc{MATLAB} implementations.

\subsection{An example of two-dimensional Helmholtz equation in homogeneous media}\label{2dsmooth}
 The exact solution (see \cite{HGA}) to the problem can be written in the closed form as $$u_{ex}(x,y)=\text{cos}(k\pi y)(A_1e^{-i\omega_x x}+A_2e^{i\omega_x
x})$$ where $\omega_x=\sqrt{\omega^2-(k\pi)^2}$, $\Omega=[0,1]\times[0,1]$, and coefficients
$A_1$ and $A_2$ satisfy the equation
\begin{equation}
\left( {\begin{array}{cc} \omega_x & -\omega_x \\
(\omega-\omega_x)e^{-2i\omega_x} & (\omega+\omega_x)e^{2i\omega_x}
\end{array} }
\right)
 \left ( {\begin{array}{c} A_1 \\ A_2
\end{array}}
\right ) = \left ( {\begin{array}{c}
-i   \\
0
\end{array}}
\right ).
\end{equation}

Then the boundary data is computed by
\be \label{2dg}
g_{ex} = (\partial_\text{\bf n} + \text{i}\omega)u_{ex}.
\en

The solution respectively represents propagating modes and
evanescent modes when the mode number $k$ is below the cut-off value
$k\leqslant k_{\text{cut-off}}={\omega\over\pi}$ and up the cut-off
value $k> k_{\text{cut-off}}.$ To be simple, we only compute the
approximate solution for the highest propagating mode with
$k=\omega/\pi-1$ in the following tests.

\begin{figure}[H]
%\vspace{-2cm}thb
\begin{center}
\begin{tabular}{cc}
\epsfxsize=0.5\textwidth\epsffile{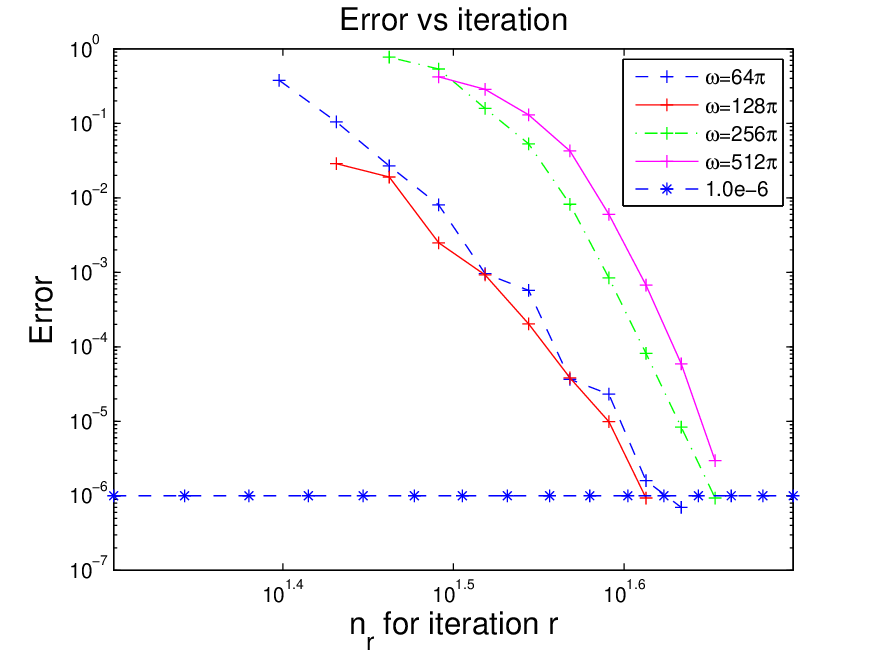}&
\epsfxsize=0.5\textwidth\epsffile{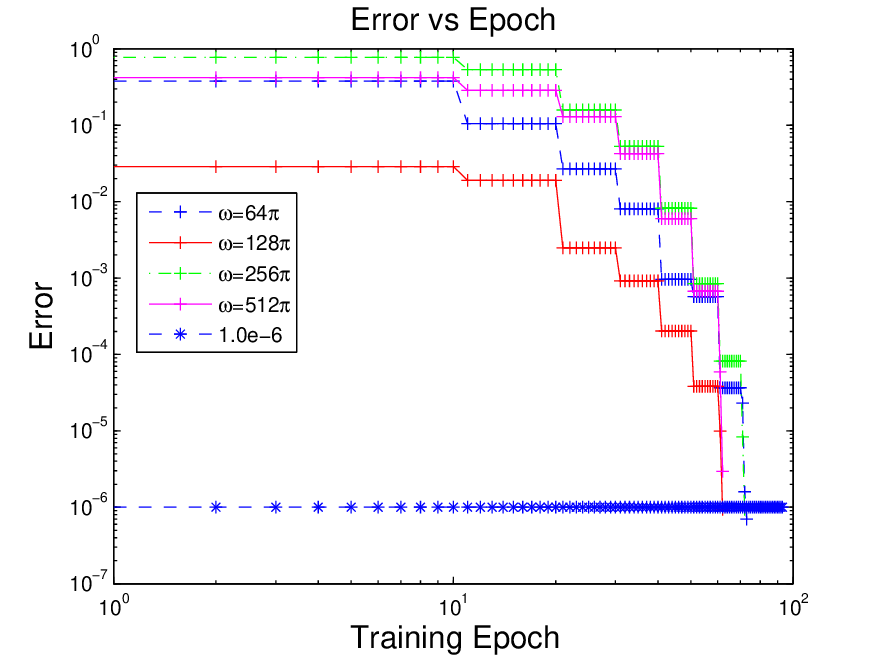}\\
\end{tabular}
\end{center}
 \caption{Displacement of a string (Helmholtz equation in two dimensions). (Left) Errors at each quasi-minimization iteration. (Right) The progress of the loss function within each quasi-minimization iteration.  }
\label{2dmul_larwave_nn1}
\end{figure}

Figure \ref{2dmul_larwave_nn1} shows the true errors $|||u-u_{r-1}|||$ at the end of each quasi-minimization iteration. We also provide the analogous results after each training epoch. Set the number of training epoches as 10 in Algorithm 4.2. The other parameters are set as follows.
\beq \nonumber
%\omega& = & 16\pi, 32\pi, ~ h=\frac{1}{8}, ~n_r = 4r+1, ~tol = 10^{-6}; \cr
\omega& = & 64\pi, ~~~~~ h=\frac{1}{16}, ~n_r = 2r+23, ~tol = 10^{-6};\cr
\omega& = & 128\pi, ~~~~~ h=\frac{1}{32},~n_r = 2r+25, ~tol = 10^{-6};\cr
\omega& = & 256\pi, ~~~~~h=\frac{1}{64},~n_r = 2r+27, ~tol = 10^{-6}; \cr
\omega& = & 512\pi, ~~~~~h=\frac{1}{128},~n_r = 2r+29, ~tol = 10^{-5}.
\eq
In general, we choose $h\approx \mathcal{O}(\frac{4\pi}{\omega})$. Here we need to use relatively great $n_r$ since we have chosen relatively coarse meshes.

%\textcolor{red}{hp-refinement strategy}
Numerical results validate that, when choosing $h\approx \mathcal{O}(\frac{4\pi}{\omega})$ and gradually increasing the width $n_r$ of the network at every iteration step of Algorithm 4.1, the resulting approximate solutions
can reach the given accuracy. The main reason for selecting parameters in this way is that the plane wave approximation error is algebraically convergent with respect to $h$ and exponentially convergent with respect to $n_r$, as described by Lemma 3.1. For relatively coarse meshes and great number $n_r$ of plane wave directions, Algorithm 4.1 exhibits surprisingly strong stability when increasing the wave number $\omega$.
In addition, it can be seen that, only no more than ten iteration counts of Algorithm 4.1 may guarantee the desired accuracy of the approximations. Besides, we observe that the
initial several steps of Algorithm 4.1 can reduce the error substantially, with further sharp decrease to a given precision.

\begin{figure}[H]
\begin{center}
\begin{tabular}{ccc}
\epsfxsize=0.3\textwidth\epsffile{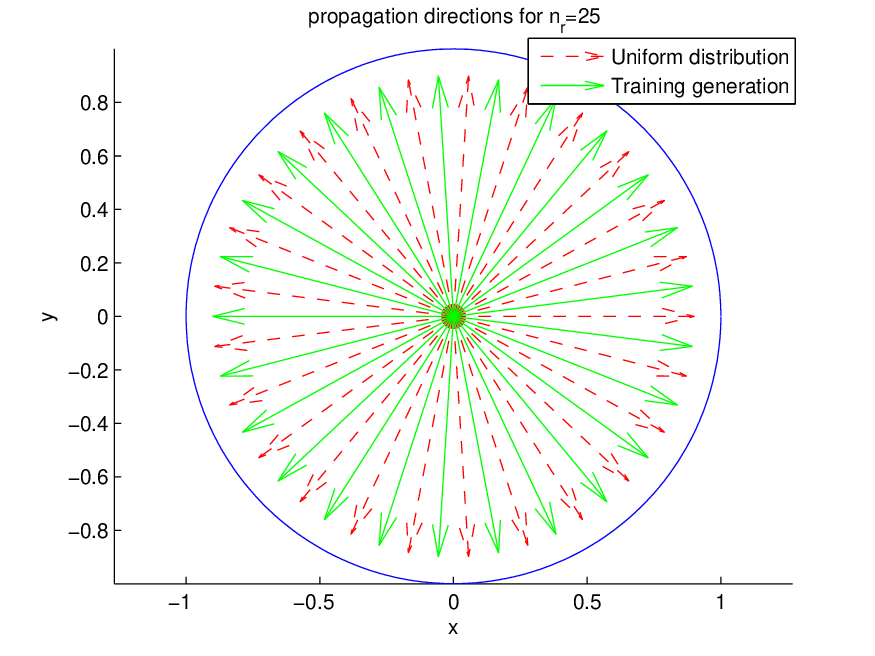}&
\epsfxsize=0.3\textwidth\epsffile{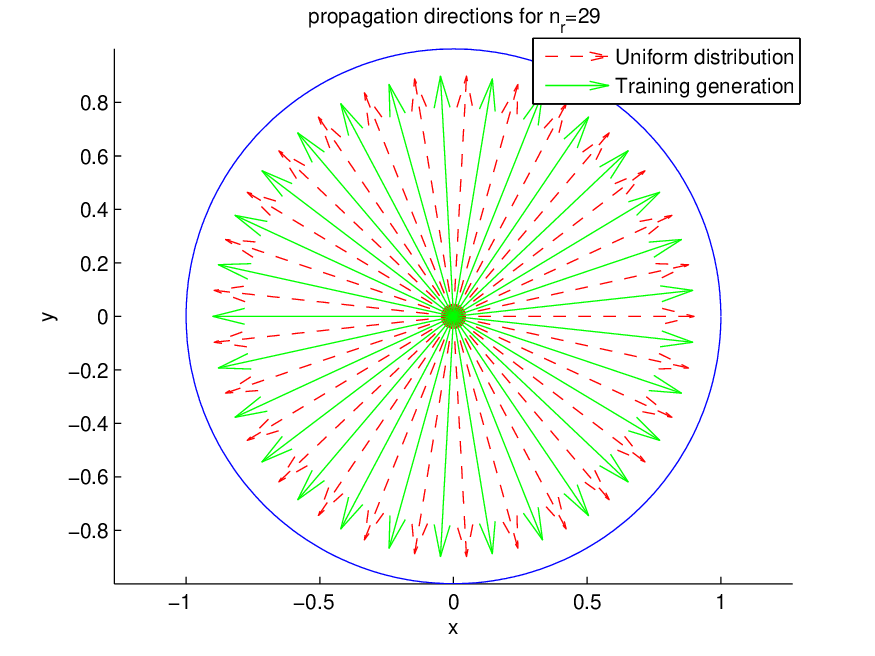}&
\epsfxsize=0.3\textwidth\epsffile{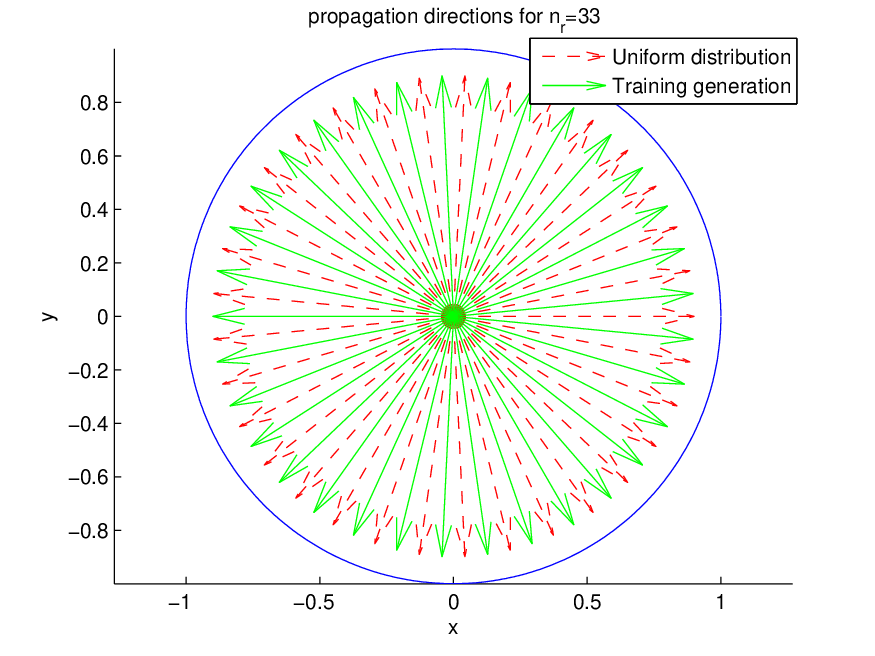}\\
\end{tabular}
\end{center}
 \caption{The trained set $\Phi^r$ in the lower-left element defining the final propagation angles of $\xi_r^{\theta}$ at several stages for $r = 1, 3, 5$. }
\label{2dhelmdirec}
\end{figure}

Figure \ref{2dhelmdirec} shows the trained set $\Phi^r$ in the lower-left element defining the final direction angles of $\xi_r^{\theta}$ at several stages of the Algorithm 1.

Compared with the uniform distribution of the direction angles employed by the PWLS method, the proposed discontinuous neural network indeed adaptively generates the desired plane wave directions.

 %\textcolor{red}{in energy-norm approximation error of $\mathcal{O}(10^{-6})$.}

Next, we would like to compare the proposed DPWNN with the PWLS method \cite{hy}. Let the number of degree of freedoms per each element in the PWLS method be the same as the width of the discontinuous network in the final
step of Algorithm 4.1. We set the maximum number of training epoches as 2 in Algorithm 4.2. Table \ref{2dhelmcomtable} shows the comparison of the errors of the resulting approximations and the computing time between two methods.

\vskip 0.1in
\begin{center}
       \tabcaption{}\vskip -0.3in
\label{2dhelmcomtable}
       Comparison of errors of approximations and the computing time with respect to $\omega$.  \vskip 0.1in
\begin{tabular}{|c|c|c|c|c|c|c|} \hline
  \multicolumn{2}{|c| } {  \(\omega\) } & $32\pi$ & $64\pi$ &  $128\pi$ & $256\pi$ & $512\pi$ \\ \hline
\multirow{2}*{ \text{DPWNN} } & \text{Error}   &  3.00e-7  & 7.00e-7   &  9.35e-7  & 9.34e-7 & 2.97e-6 \\
& \text{Time}   &  8.09e+1  &  9.12e+2  &  2.25e+3 &  5.11e+3 &  1.67e+4 \\ \hline
 \multirow{2}*{ \text{PWLS} }  & \text{Error}   &  1.72e-6 &   7.73e-6 &  4.84e-6  & 6.01e-6  &  1.89e-5 \\
& \text{Time}   & 3.96e+1   &  4.32e+2   &  9.57e+2  & 2.67e+3 & 9.25e+3   \\ \hline
   \end{tabular}
     \end{center}

\vskip 0.1in

\begin{figure}[H]
%\vspace{-2cm}thb
\begin{center}
\begin{tabular}{cc}
\epsfxsize=0.5\textwidth\epsffile{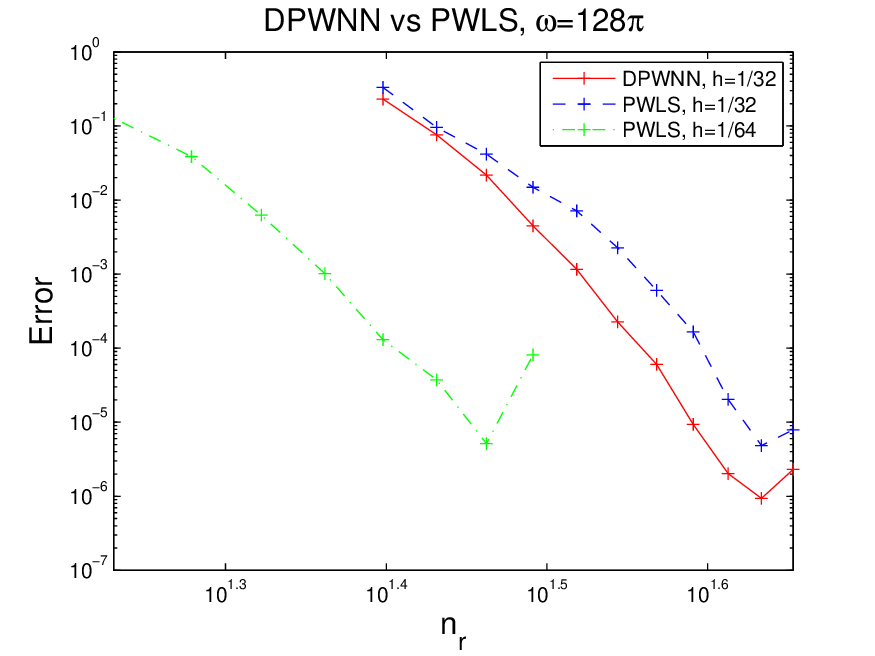}&
\epsfxsize=0.5\textwidth\epsffile{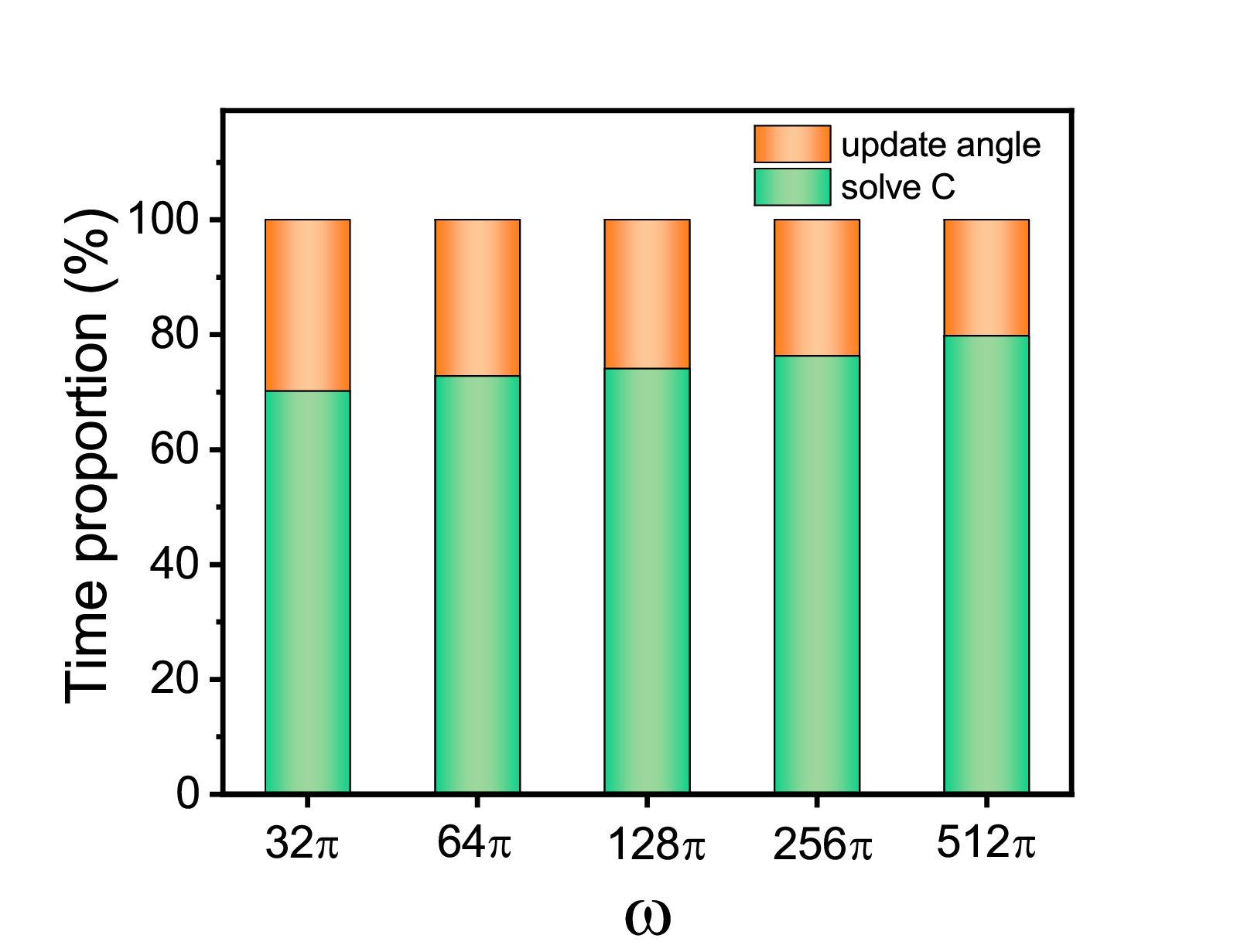}\\
\end{tabular}
\end{center}
 \caption{Comparison between the DPWNN and the PWLS method (Helmholtz equation in two dimensions). (Left) Errors at each quasi-minimization iteration. (Right) Proportion of computing time w.r.t. $\omega$.  }
\label{2dhelm_refine_compare}
\end{figure}

Figure \ref{2dhelm_refine_compare} (Left) shows the comparison of the errors of the resulting approximations with respect to $n_r$ for the case of a fixed $\omega=128\pi$. Figure \ref{2dhelm_refine_compare} (Right) shows the proportion of computing time spent on solving the linear systems (\ref{galerlsq}) and updating propagation directions, respectively.

The results listed in Table \ref{2dhelmcomtable} indicate that the approximation generated by the DPWNN method is much more accurate than that generated by the original PWLS method. This is mainly caused by adaptive
superpositions with corrected plane wave propagation directions generated by iterative algorithms,  as shown in Figure \ref{2dhelmdirec}. The results listed in Table \ref{2dhelmcomtable} also indicate that
the computational cost for the entire iteration of DPWNN is almost twice that of the PWLS method. This is mainly because the maximum number of training epoches is set as 2 in Algorithm 4.2, and the computational expense in Algorithm 4.1 mainly spends on the calculation of the entire sequence of activation coefficients $c^{r}$ defined by the linear system (\ref{galerlsq}), which can be illustrated by Figure \ref{2dhelm_refine_compare} (Right).

A natural question is whether the accuracies of the approximations generated by the PWLS method can be essentially improved if further increasing the cost, for example, decreasing the meshsize $h$? Unfortunately, the answer is negative: the approximations generated by the PWLS method can not reach the high accuracies of the approximate solutions generated by the DPWNN method since the standard plane wave approaches suffer from serious ill-conditioning and are heavily unstable \cite{ref11,ref21,HMPsur,HGA}, see Figure \ref{2dhelm_refine_compare} (Left).
%thus without an appropriate preconditioning it is impossible to obtain meaningful results for sufficiently large $n_r$. We refer to \cite{hut,HMPsur} for a further remedy of this issue.

\subsection{An example of three-dimensional Helmholtz equation}
The following test problem consists of a point source and associated
boundary conditions for homogeneous Helmholtz equations (see
\cite{HKM}):
\begin{equation}
\begin{split}
&u(r,r_0)= {1\over 4\pi} {e^{i\omega|r-r_0|}\over |r-r_0|} ~~ \text{in} ~~\Omega,  \\
&{\partial u\over \partial {\bf n}}+i\omega u=g \quad
\text{over}\quad\partial\Omega, \\
\end{split}\label{fig1}
\end{equation}
in a cubic computational domain $\Omega=[0,1]\times[0,1]\times[0,1]$. The location of the source is off-centred at
$r_0=(-1,-1,-1)$ and $r=(x,y,z)$ is an observation point. Note that the
analytic solution of the Helmholtz equation with such boundary
condition has a singularity at $r=r_0$.

\begin{figure}[H]
%\vspace{-2cm}thb
\begin{center}
\begin{tabular}{cc}
\epsfxsize=0.5\textwidth\epsffile{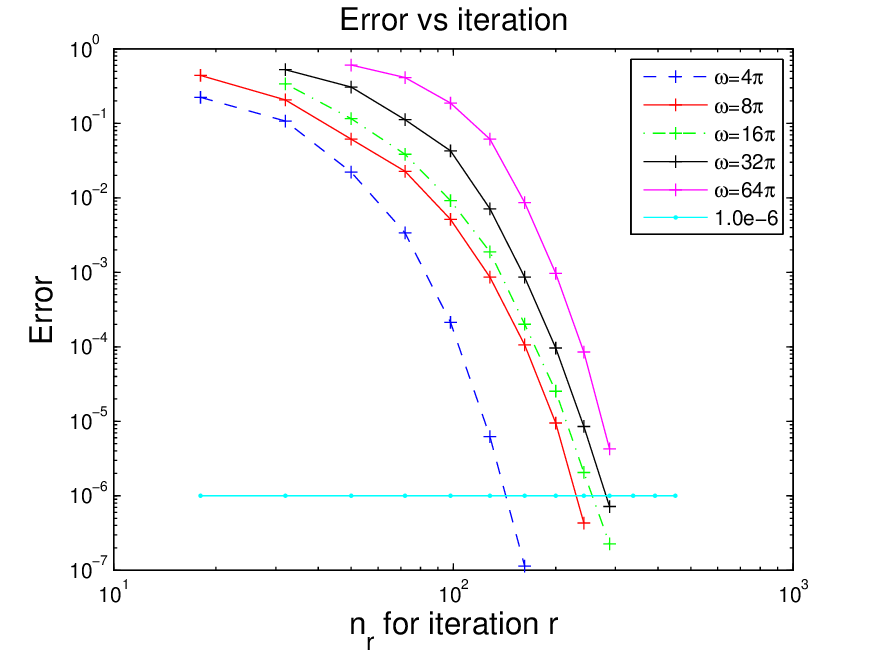}&
\epsfxsize=0.5\textwidth\epsffile{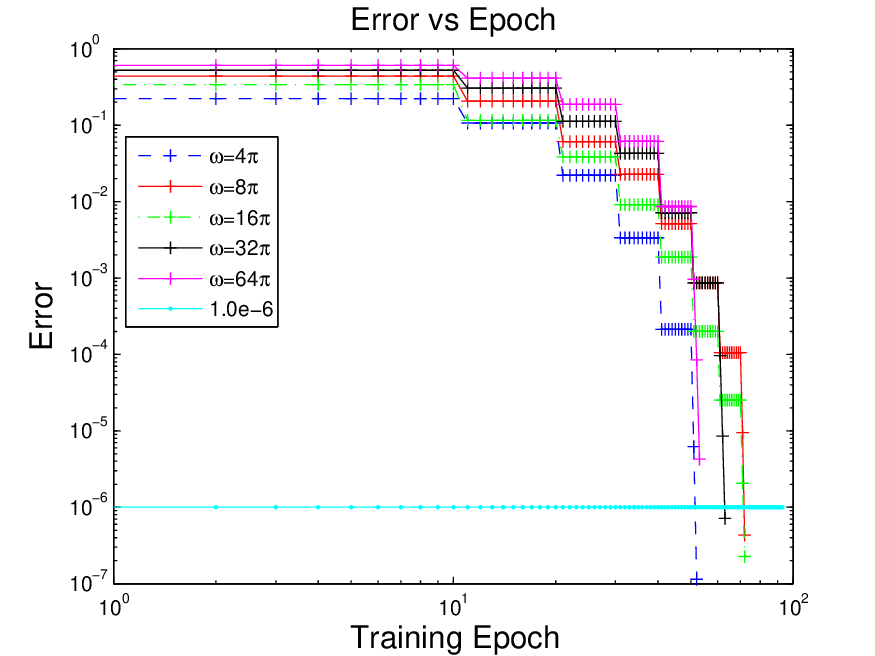}\\
\end{tabular}
\end{center}
 \caption{Displacement of a string (Helmholtz equation in three dimensions). (Left) Errors at each quasi-minimization iteration. (Right) The progress of the loss function within each quasi-minimization iteration.  }
\label{3dhelm_nn}
\end{figure}

Figure \ref{3dhelm_nn} shows the true errors $|||u-u_{r-1}|||$ at the end of each quasi-minimization iteration. We also provide the analogous results after each training epoch. Set the number of training epoches as 10 in Algorithm 4.2. The other parameters are set as follows.
\beq \nonumber
\omega& = & 4\pi, 8\pi, ~ h=\frac{1}{2}, ~m^{\ast}_{r} = r+2, ~tol = 10^{-6}; \cr
\omega& = & 16\pi, ~~~~~ h=\frac{1}{4}, ~m^{\ast}_{r} = r+3, ~tol = 10^{-6};\cr
\omega& = & 32\pi, ~~~~~ h=\frac{1}{8},~m^{\ast}_{r} = r+3, ~tol = 10^{-6};\cr
\omega& = & 64\pi, ~~~~~ h=\frac{1}{16},~m^{\ast}_{r} = r+4, ~tol = 10^{-5}.
\eq
Generally, we choose $h\approx \mathcal{O}(\frac{4\pi}{\omega})$.

Numerical results validate that, when choosing $h\approx \mathcal{O}(\frac{4\pi}{\omega})$ and gradually increasing the width $n_r$ of the network at every iteration step of Algorithm 4.1, the resulting approximate solutions
can reach the given accuracy. As in the last example, Algorithm 4.1 exhibits strong stability when increasing the wave number $\omega$, and
only no more than ten iteration counts of Algorithm 4.1 may guarantee the desired accuracy of approximations.

%For relatively coarse meshes and great number $n_r$ of plane wave directions, Algorithm 4.1 exhibits surprisingly strong stability when increasing the wave number $\omega$.
%In addition, it can be seen that, only no more than ten iteration counts of Algorithm 4.1 may guarantee the desired accuracy of approximations.

\begin{figure}[H]
%\vspace{-2cm}thb
\begin{center}
\begin{tabular}{cc}
\epsfxsize=0.4\textwidth\epsffile{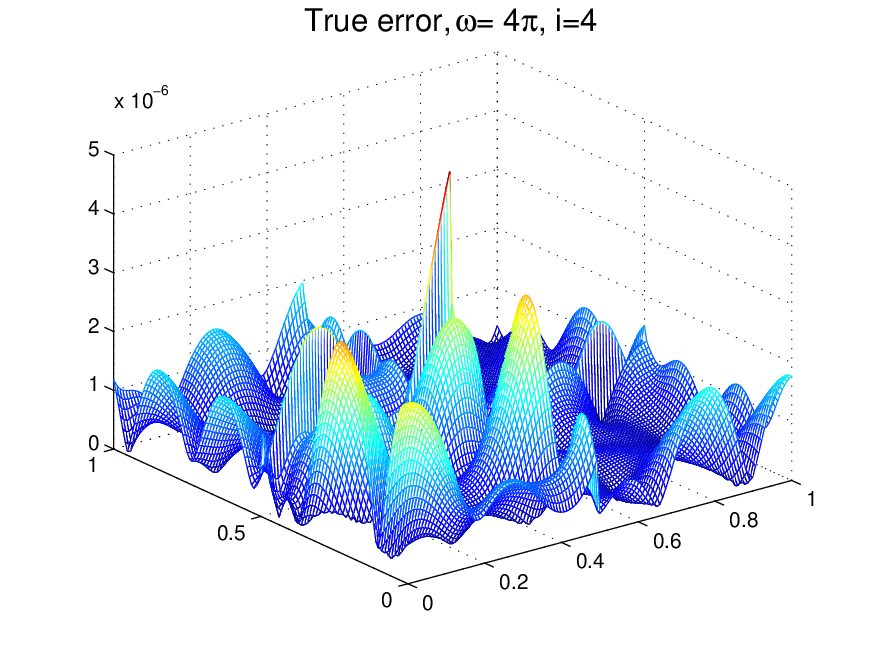}&
\epsfxsize=0.4\textwidth\epsffile{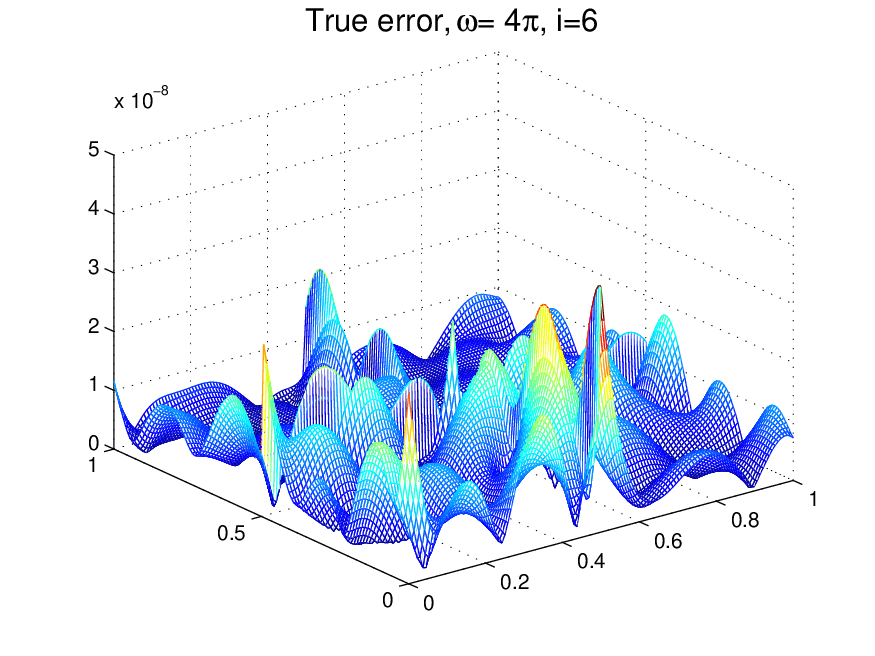}\\
\end{tabular}
\end{center}
 \caption{Displacement of a membrane (Helmholtz equation in three dimensions). Exact error $|u-u_{i-1}|$ for i = 4, 6. }
\label{3dhenn2}
\end{figure}

Figure \ref{3dhenn2} shows the exact error $|u-u_{r-1}|$ at several stages of the algorithm for the case of cross-section $z=0.5$.

It is apparent that initially, the low frequency components of the error are learned, with later iterations learning the high frequency error components, and the pointwise error is consistent with the quasi-minimization iteration error in Figure \ref{3dhelm_nn}.

\begin{figure}[H]
%\vspace{-2cm}thb
\begin{center}
\begin{tabular}{cc}
\epsfxsize=0.5\textwidth\epsffile{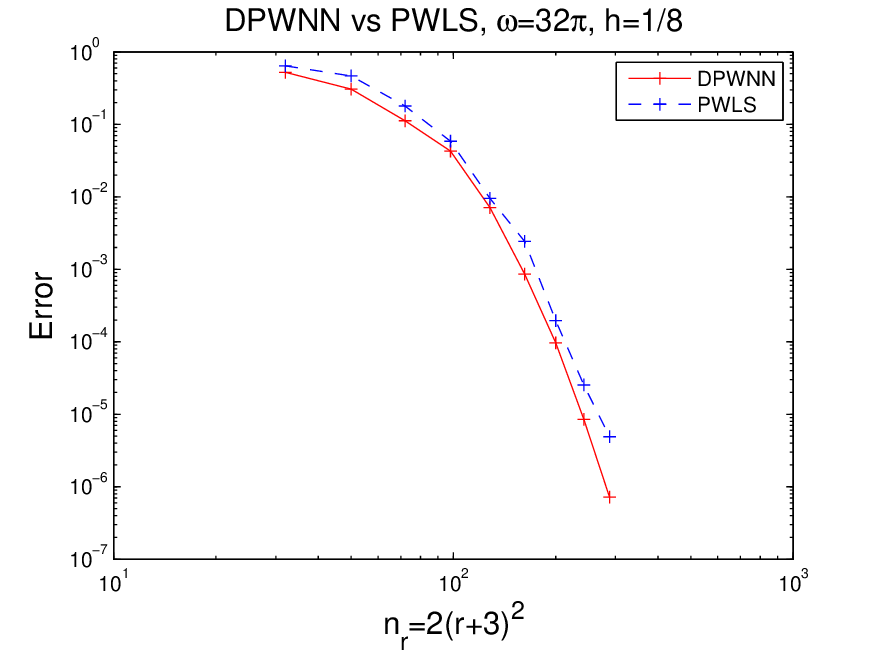}&
\epsfxsize=0.5\textwidth\epsffile{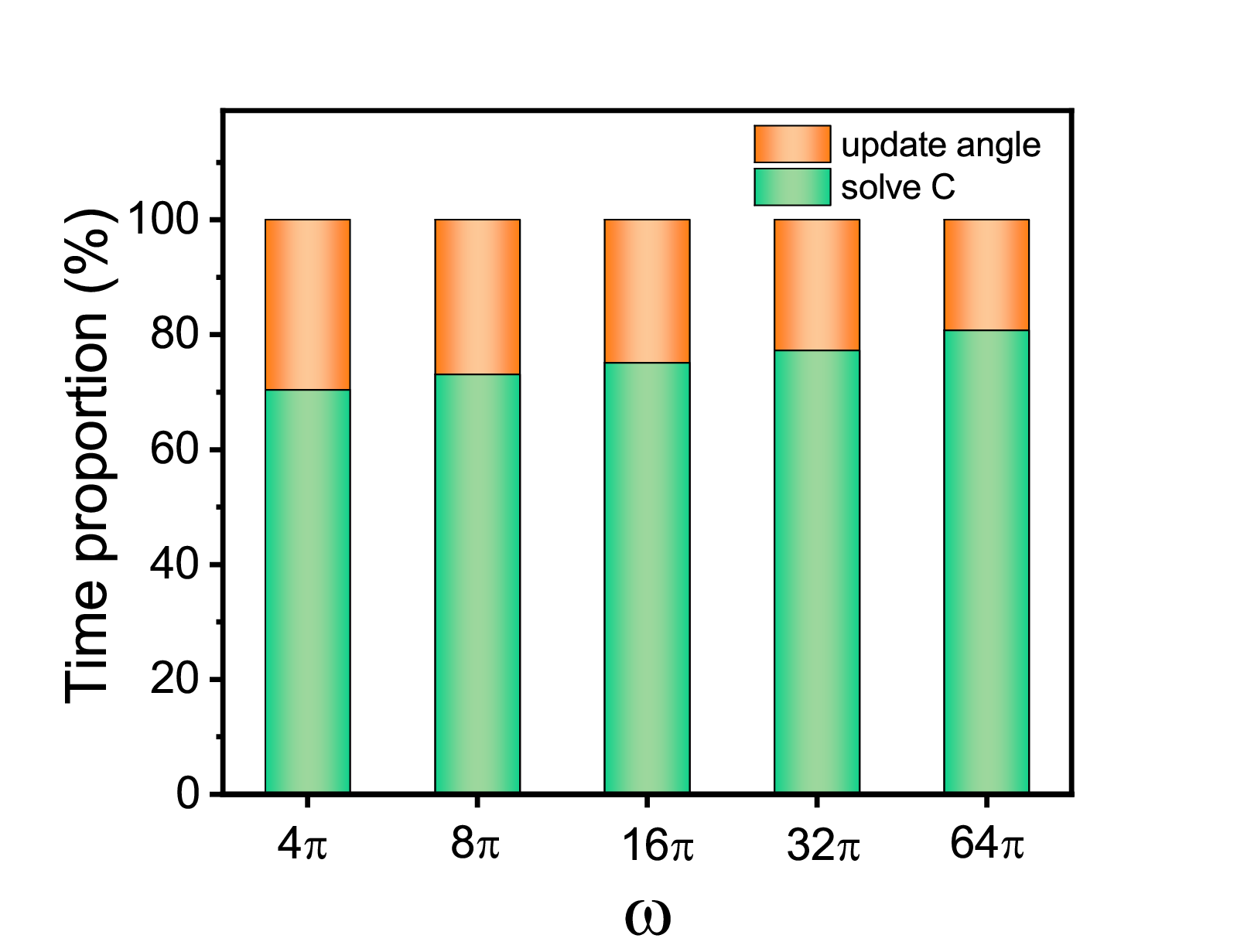}\\
\end{tabular}
\end{center}
 \caption{Comparison between the DPWNN and the PWLS method (Helmholtz equation in three dimensions). (Left) Errors at each quasi-minimization iteration. (Right) Proportion of computing time w.r.t. $\omega$.  }
\label{3dhelm_refine_compare}
\end{figure}

Let the number of degree of freedoms per each element in the PWLS method be the same as the width of the discontinuous network in the final step of Algorithm 4.1. We set the maximum number of training epoches as 2 in Algorithm 4.2. Table \ref{3dhelmcomtable} shows the comparison of the errors of the resulting approximations and the computing time between two methods. Figure \ref{3dhelm_refine_compare} (Left) shows the comparison of the errors of the resulting approximations with respect to $n_r$ for the case of a fixed $\omega=32\pi$. Figure \ref{3dhelm_refine_compare} (Right) shows the proportion of computing time spent on solving the linear systems (\ref{galerlsq}) and updating propagation directions, respectively.

\vskip 0.1in
\begin{center}
       \tabcaption{}\vskip -0.3in
\label{3dhelmcomtable}
       Comparison of errors of approximations and the computing time with respect to $\omega$.  \vskip 0.1in
\begin{tabular}{|c|c|c|c|c|c|c|} \hline
  \multicolumn{2}{|c| } {  \(\omega\) } & $4\pi$ & $8\pi$ &  $16\pi$ & $32\pi$ & $64\pi$ \\ \hline
\multirow{2}*{ \text{DPWNN} } & \text{Error}   &  1.14e-7  & 4.32e-7   &  2.28e-7  & 7.18e-7 & 4.27e-6 \\
& \text{Time}   & 2.35e+2   &  4.54e+2  &  5.48e+3  & 1.38e+4  & 1.14e+5 \\ \hline
 \multirow{2}*{ \text{PWLS} }  & \text{Error}     &  6.29e-7   & 3.02e-6    &  1.53e-6  & 4.90e-6  &  2.96e-5   \\
& \text{Time} & 1.23e+2 & 2.33e+2  & 2.70e+3  &  6.51e+3  & 5.68e+4  \\ \hline
   \end{tabular}
     \end{center}

\vskip 0.1in

The results listed in Table \ref{3dhelmcomtable} also indicate that the approximation generated by the DPWNN method is much more accurate than that generated by the original PWLS method,
%owing to adaptive superpositions with corrected plane wave propagation directions generated by iterative algorithms,
as shown in Figure \ref{3dhelm_refine_compare} (Left).

%Moreover, from Table \ref{3dhelmcomtable}, the computational cost for the entire iteration of DPWNN is almost twice that of the PWLS method. This is mainly because the maximum number of training epoches is set as 2 in Algorithm %4.2, and the main computational expense in Algorithm 4.1 mainly involves solving the entire sequence of activation coefficients $c^{r}$ defined by the linear system (\ref{galerlsq}), which can be illustrated by Figure %\ref{3dhelm_refine_compare} (Right).

\subsection{An example of three-dimensional Maxwell's equations in homogeneous media}\label{3dhomomax}

We compute the electric field due to an electric dipole source at
the point \({\bf x}_0=(0.6,0.6,0.6)\). The dipole point source can
be defined as the solution of a homogeneous Maxwell system
(\ref{maxeq}). The exact solution of the problems is
\begin{equation}
{\bf E}_{\text{ex}}=-\text{i}\omega I \phi({\bf x},{\bf x}_0){ \bf a}
+\frac{I}{{ i}\omega\varepsilon} \nabla(\nabla\phi\cdot{ \bf a}),
\end{equation}
where
$$ \phi({\bf x},{\bf x}_0) = \frac{\text{exp}(\text{i}\omega\sqrt{\varepsilon}|{\bf x}-{\bf x}_0|)}{4\pi|{\bf x}-{\bf x}_0|}  $$
and \(\Omega=[-0.5,0.5]^3\). Then the boundary data is computed by
\be \label{3dg}
{\bf g}_{ex} =
-{\bf E}_{\text{ex}}\times {\bf n}+\frac{\sigma}{i\omega\mu}((\nabla\times {\bf E}_{\text{ex}})\times  {\bf n})\times {\bf n}.
\en

To keep the exact solution smooth in
\(\Omega\), we move the singularity \( {\bf x}_0\) from
\((0.2,0.2,0.2)\) (see \cite{hmm}) in the computational domain
to \((0.6,0.6,0.6)\) outside the region.

\begin{figure}[H]
%\vspace{-2cm}thb
\begin{center}
\begin{tabular}{cc}
\epsfxsize=0.5\textwidth\epsffile{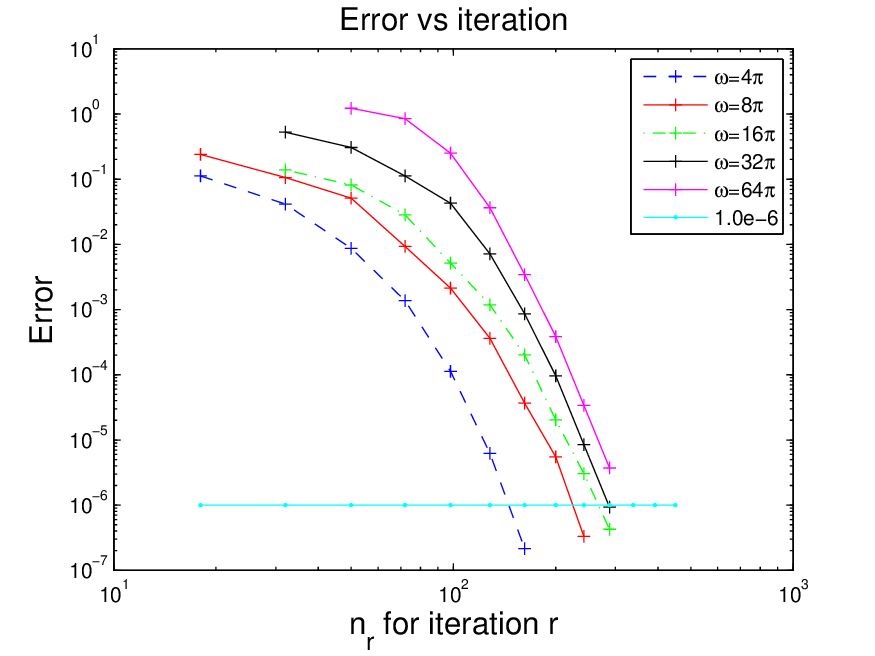}&
\epsfxsize=0.5\textwidth\epsffile{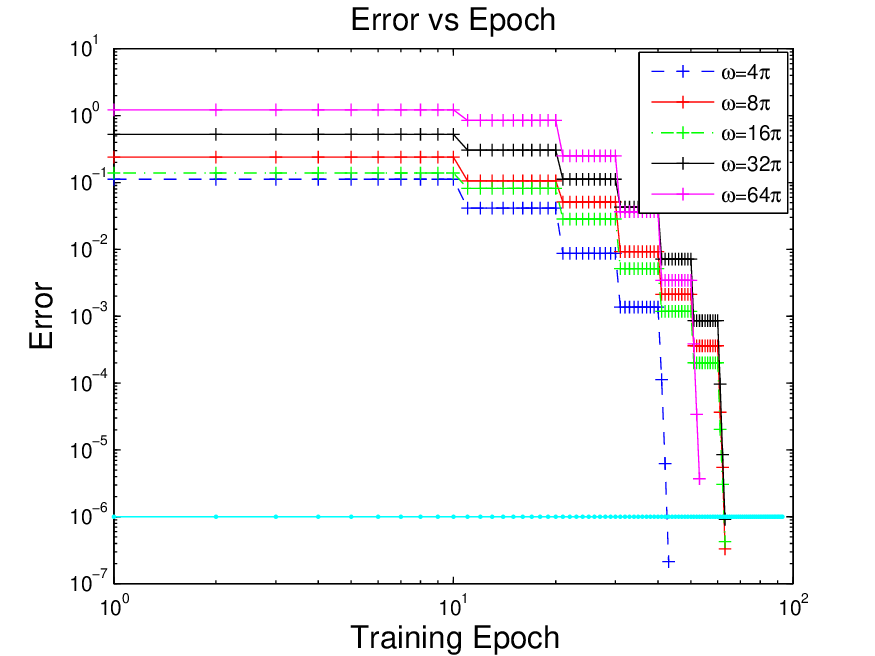}\\
\end{tabular}
\end{center}
 \caption{Displacement of a string (Maxwell's equations in three dimensions). (Left) Errors at each quasi-minimization iteration. (Right) The progress of the loss function within each quasi-minimization iteration.  }
\label{3dmax_nn}
\end{figure}

Figure \ref{3dmax_nn} shows the true errors $|||{\bf E}-{\bf E}_{r-1}|||$ at the end of each quasi-minimization iteration. We also provide the analogous results after each training epoch. Set the number of training epoches as 10 in Algorithm 4.2. The other parameters are set as follows.
\beq \nonumber
\omega& = & 4\pi, 8\pi, ~ h=\frac{1}{2}, ~m^{\ast}_{r} = r+2, ~tol = 10^{-6}; \cr
\omega& = & 16\pi, ~~~~~ h=\frac{1}{4}, ~m^{\ast}_{r} = r+3, ~tol = 10^{-6};\cr
\omega& = & 32\pi, ~~~~~ h=\frac{1}{8},~m^{\ast}_{r} = r+3, ~tol = 10^{-6};\cr
\omega& = & 64\pi, ~~~~~ h=\frac{1}{16},~m^{\ast}_{r} = r+4, ~tol = 10^{-5}.
\eq
Generally, we choose $h\approx \mathcal{O}(\frac{4\pi}{\omega})$.

\iffalse;\cr
\omega& = & 128\pi, ~~~~~ h=\frac{1}{32},~m^{\ast}_{r} = r+4, ~tol = 10^{-5};\cr
\omega& = & 256\pi, ~~~~~ h=\frac{1}{64},~m^{\ast}_{r} = r+5, ~tol = 10^{-4};\cr
\omega& = & 512\pi, ~~~~~ h=\frac{1}{128},~m^{\ast}_{r} = r+5, ~tol = 10^{-4}.
\fi

Numerical results validate that, when choosing $h\approx \mathcal{O}(\frac{4\pi}{\omega})$ and gradually increasing the width $n_r$ of the network at every iteration step of Algorithm 4.1, the resulting approximate solutions
can reach the given accuracy. As in the example tested in subsection 6.1, Algorithm 4.1 exhibits strong stability when increasing the wave number $\omega$, and only no more than ten iteration counts of Algorithm 4.1 may guarantee the desired accuracy of the approximations.

\begin{figure}[H]
\begin{center}
\begin{tabular}{ccc}
\epsfxsize=0.3\textwidth\epsffile{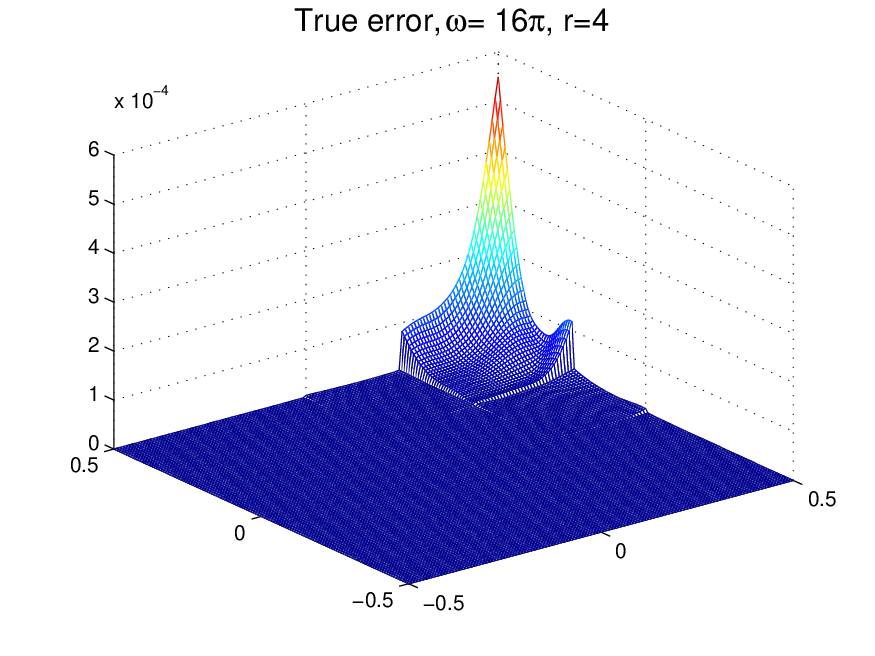}&
\epsfxsize=0.3\textwidth\epsffile{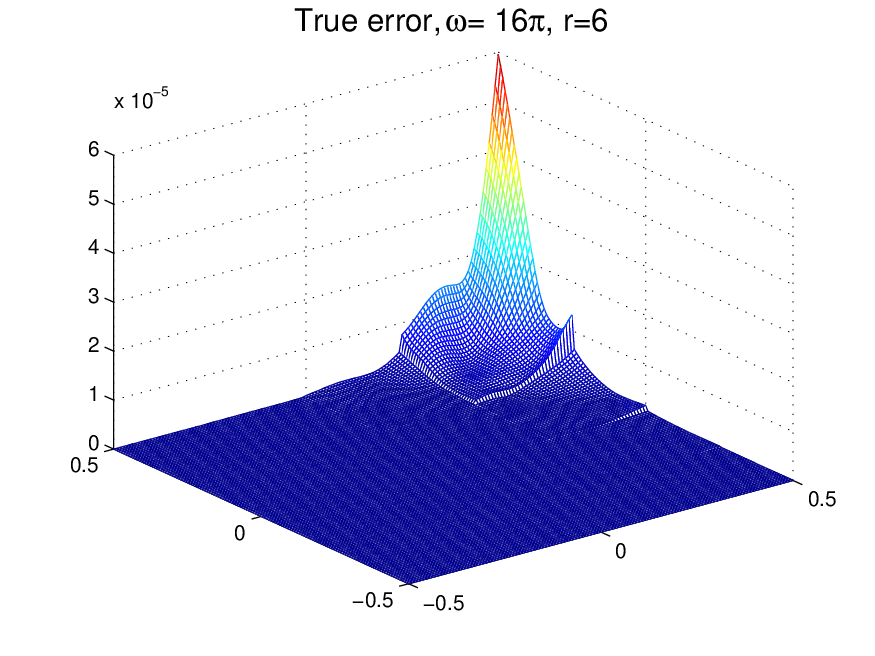}&
\epsfxsize=0.3\textwidth\epsffile{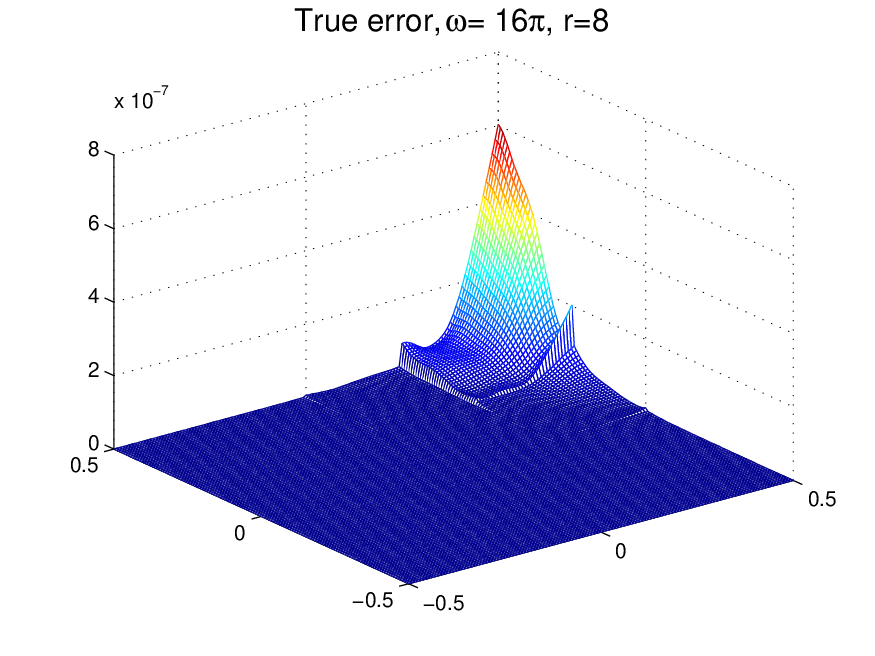}\\
\end{tabular}
\end{center}
 \caption{Displacement of a membrane (Maxwell's equations in three dimensions). Exact error $|{\bf E}-{\bf E}_{r}|$ for r = 4, 6, 8. }
\label{3dmax2}
\end{figure}

Figure \ref{3dmax2} shows the exact error $|{\bf E}-{\bf E}_{r}|$ at several stages of the algorithm for the case of cross-section $z=0$.

It is apparent that the pointwise error is consistent with the quasi-minimization iteration error in Figure \ref{3dmax_nn}. Besides, we can find that, the closer to the singularity point ${x}_0=(0.2,0.2,0.2)$, the larger the pointwise error. To improve this situation, the local refinement on the meshwidth $h$ and the degree of freedoms $n_r$ per each element can be employed in the future.

\begin{figure}[H]
%\vspace{-2cm}thb
\begin{center}
\begin{tabular}{cc}
\epsfxsize=0.5\textwidth\epsffile{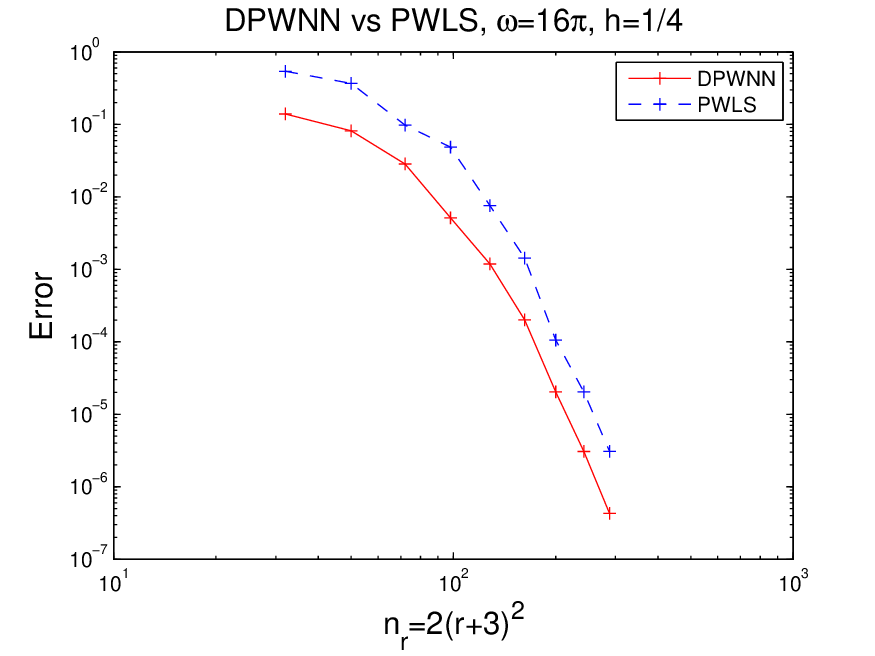}&
\epsfxsize=0.5\textwidth\epsffile{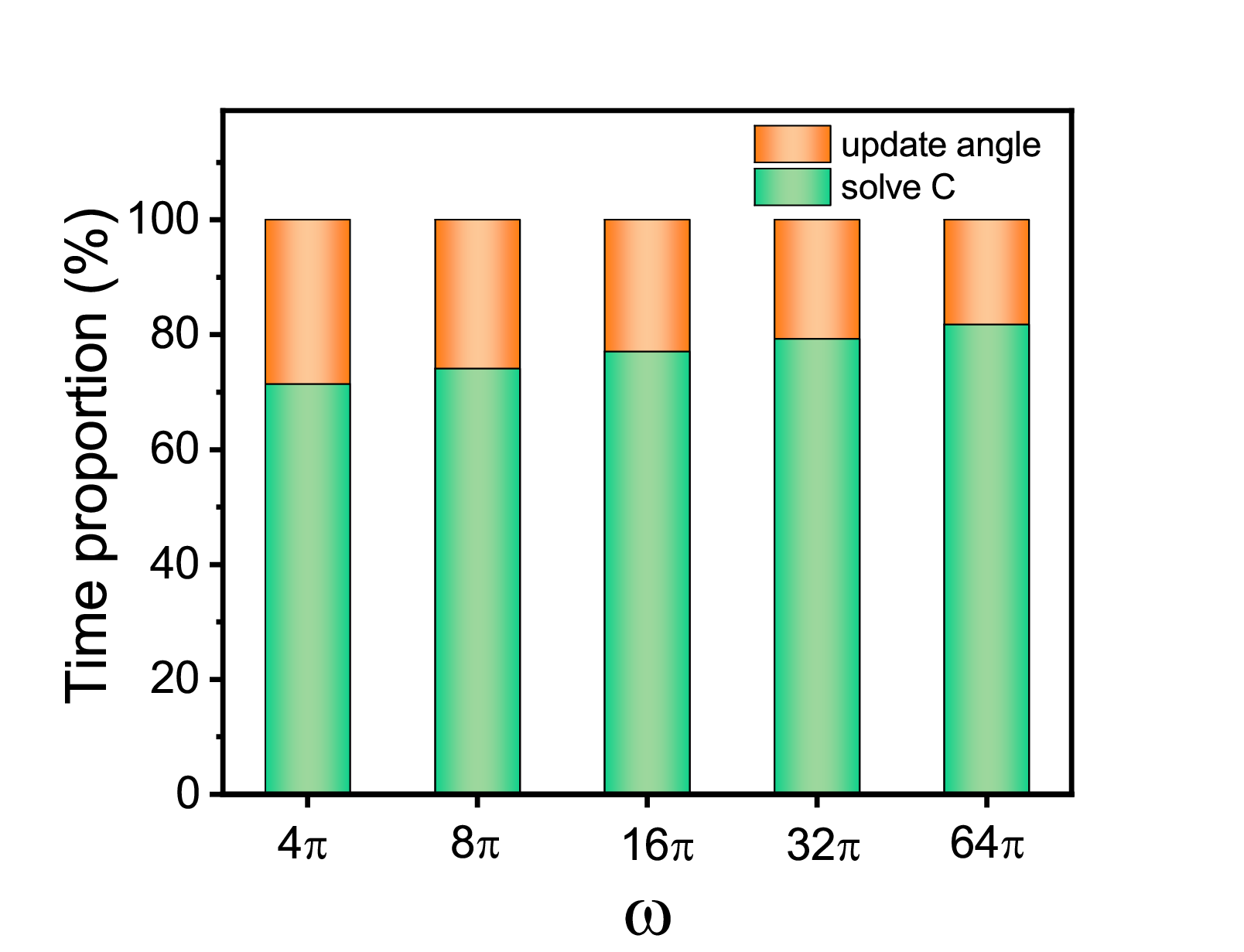}\\
\end{tabular}
\end{center}
 \caption{Comparison between the DPWNN and the PWLS method (Maxwell's equation in three dimensions). (Left) Errors at each quasi-minimization iteration. (Right) Proportion of computing time w.r.t. $\omega$.  }
\label{3dmax_refine_compare}
\end{figure}

Let the number of degree of freedoms per each element in the PWLS method be the same as the width of the discontinuous network in the final step of Algorithm 4.1. We set the maximum number of training epoches as 2 in Algorithm 4.2. Table \ref{3dmaxcomtable} shows the comparison of the errors of the resulting approximations and the computing time between two methods. Figure \ref{3dmax_refine_compare} (Left) shows the comparison of the errors of the resulting approximations with respect to $n_r$ for the case of a fixed $\omega=16\pi$. Figure \ref{3dmax_refine_compare} (Right) shows the proportion of computing time spent on solving the linear systems (\ref{maxpwlsvar}) and updating propagation directions, respectively.

\vskip 0.1in
\begin{center}
       \tabcaption{}\vskip -0.3in
\label{3dmaxcomtable}
       Comparison of errors of approximations and the computing time with respect to $\omega$.  \vskip 0.1in
\begin{tabular}{|c|c|c|c|c|c|c|} \hline
  \multicolumn{2}{|c| } {  \(\omega\) } & $4\pi$ & $8\pi$ &  $16\pi$ & $32\pi$ & $64\pi$ \\ \hline
\multirow{2}*{ \text{DPWNN} } & \text{Error}   &  2.14e-7  & 3.32e-7   &  4.27e-7  & 9.18e-7 & 3.83e-6 \\
& \text{Time}   & 2.80e+2   &  5.23e+2  & 5.90e+3  & 1.50e+4  & 1.24e+5 \\ \hline
 \multirow{2}*{ \text{PWLS} }  & \text{Error}   & 1.38e-6 &   2.29e-6 &  3.08e-6  & 6.10e-6  &  2.81e-5 \\
& \text{Time}   &  1.37e+2  & 2.54e+2   & 2.81e+3  & 7.16e+3  &  5.79e+4  \\ \hline
   \end{tabular}
     \end{center}

\vskip 0.1in

The results listed in Table \ref{3dmaxcomtable} also indicate that the approximation generated by the DPWNN method is much more accurate than that generated by the original PWLS method,
%owing to adaptive superpositions with corrected plane wave propagation directions generated by iterative algorithms,
as shown in Figure \ref{3dmax_refine_compare} (Left).

% Moreover, from Table \ref{3dmaxcomtable}, the computational cost for the entire iteration of DPWNN is almost twice that of the PWLS method. This is mainly because the maximum number of training epoches is set as 2 in Algorithm %4.2, and the main computational expense in Algorithm 4.1 mainly involves solving the entire sequence of activation coefficients $c^{r}$ defined by the linear system (\ref{maxpwlsvar}), which can be illustrated by Figure %\ref{3dmax_refine_compare} (Right).

\subsection{An example of two-dimensional Helmholtz equation in inhomogeneous media}
In this subsection, we test an example in the inhomogeneous media. For the case of inhomogeneous media (i.e., $\omega$ is not a constant), it is difficult to give an analytic solution of the Helmholtz system (\ref{helm1}). As usual we replace the analytic solution by the numerical solution generated by the proposed method on the finer mesh and the larger width of the network.

%As usual we replace the analytic solution by its ``good" approximation generated by the standard finite element method with very fine grids in order to compute accuracies of the approximations generated by the proposed method.

Let the boundary data $g$ in the equation (\ref{helm1}) be chosen as the same function $g_{ex}$ given in section \ref{2dsmooth}. The other quantities are kept the same. But we choose $\omega=\omega_r$ for the right domain $x>\frac{1}{2}$ and $\omega=2\omega_r$ for the subregion $x<\frac{1}{2}$. Besides, the relaxation parameter $\alpha$ is set to be $\alpha|_{\Gamma_{kj}}=\omega_k\omega_j$, where $\omega_k=\omega|_{\Omega_k}$.

%Since the parameter $\omega$ is different from the one chosen in section \ref{2dsmooth}, the analytic function $u_{ex}$ in section \ref{2dsmooth} is not the analytic solution of the current equations.

\begin{figure}[H]
%\vspace{-2cm}thb
\begin{center}
\begin{tabular}{cc}
\epsfxsize=0.5\textwidth\epsffile{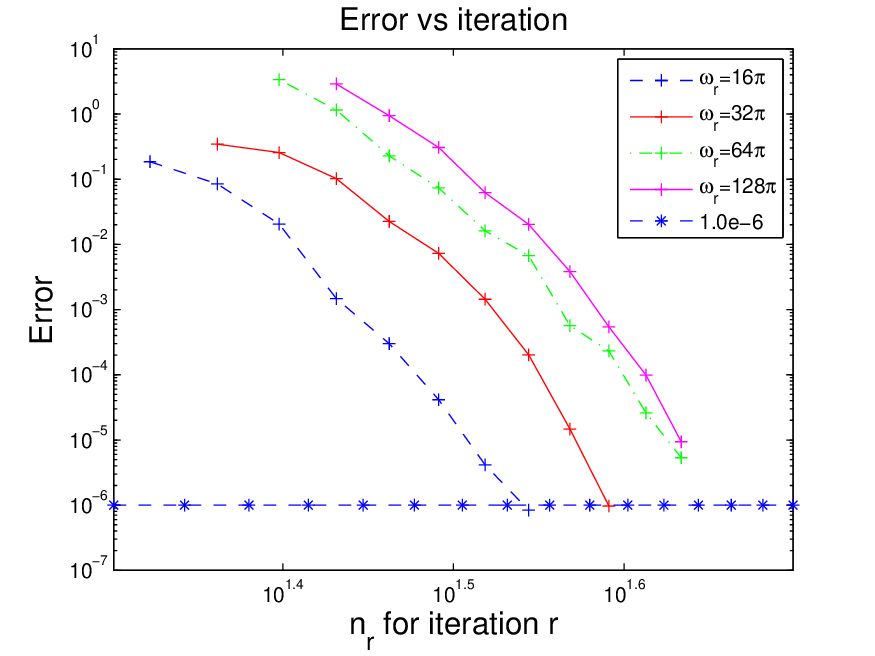}&
\epsfxsize=0.5\textwidth\epsffile{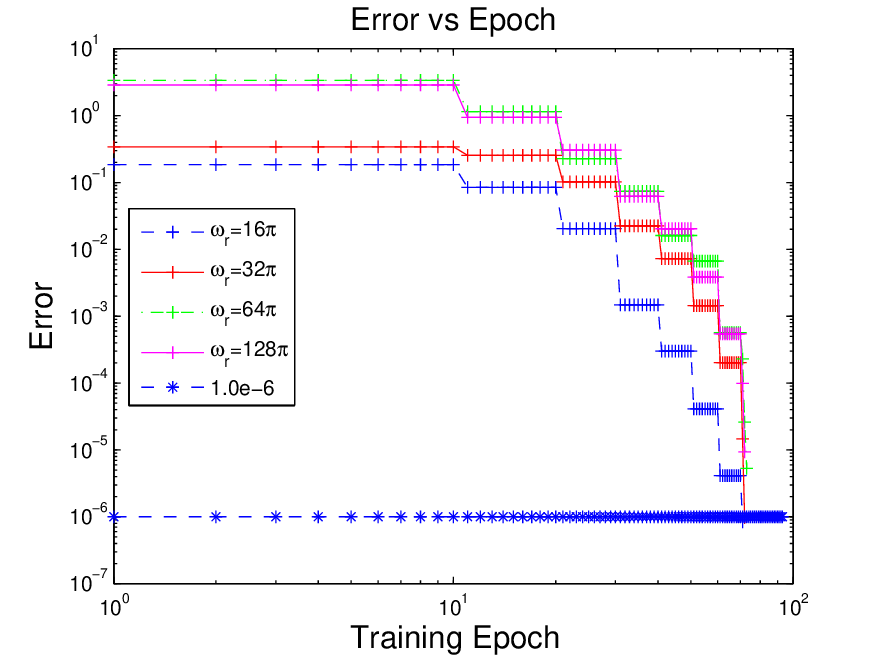}\\
\end{tabular}
\end{center}
 \caption{Displacement of a string (Helmholtz equation in two dimensions). (Left) Errors at each quasi-minimization iteration. (Right) The progress of the loss function within each quasi-minimization iteration.  }
\label{2dinhomo_mul_larwave_nn1}
\end{figure}

Figure \ref{2dinhomo_mul_larwave_nn1} shows the estimated errors of the true errors $|||u-u_{r-1}|||$.  %at the end of each quasi-minimization iteration. We also provide the analogous results after each training epoch.
The relative parameters are set as follows.
\beq \nonumber
\omega_r& = & 16\pi, ~~~~~ h=\frac{1}{4}, ~n_r = 2r+19, ~tol = 10^{-6}; \cr
\omega_r& = & 32\pi, ~~~~~ h=\frac{1}{8}, ~n_r = 2r+21, ~tol = 10^{-6}; \cr
\omega_r& = & 64\pi, ~~~~~ h=\frac{1}{16}, ~n_r = 2r+23, ~tol = 10^{-5};\cr
\omega_r& = & 128\pi, ~~~~~ h=\frac{1}{32},~n_r = 2r+25, ~tol = 10^{-5}.
%\omega& = & 256\pi, ~~~~~h=\frac{1}{64},~n_r = 2r+27, ~tol = 10^{-3}.
\eq
%In general, we choose $h\approx \mathcal{O}(\frac{4\pi}{\omega})$. Here we need to use relatively great $n_r$ since we have chosen relatively coarse meshes.

Numerical results validate that, the resulting approximate solutions
can reach the given accuracy. %when choosing $h\approx \mathcal{O}(\frac{4\pi}{\omega})$ and gradually increasing the width $n_r$ of the network at every iteration step of Algorithm 4.1,

%Next, we would like to compare the proposed DPWNN with the PWLS method \cite{hy}. %Let the number of degree of freedoms per each element in the PWLS method be the same as the width of the discontinuous network in the final step of Algorithm 4.1. We set the maximum number of training epoches as 2 in Algorithm 4.2.
  Table \ref{2dinhelmcomtable} shows the comparison of the errors of the resulting approximations and the computing time between the proposed DPWNN and the PWLS method.

\vskip 0.1in
\begin{center}
       \tabcaption{}\vskip -0.3in
\label{2dinhelmcomtable}
       Comparison of errors of approximations and the computing time with respect to $\omega_r$.  \vskip 0.1in
\begin{tabular}{|c|c|c|c|c|c|} \hline
  \multicolumn{2}{|c| } {  \(\omega_r\) } & $16\pi$ & $32\pi$ & $64\pi$ &  $128\pi$  \\ \hline
\multirow{2}*{ \text{DPWNN} } & \text{Error}   & 8.35e-7  & 9.70e-7  & 5.32e-6  &  9.35e-6   \\
& \text{Time}   & 3.02e+1  & 8.76e+1  &  9.36e+2  &  2.55e+3 \\ \hline
 \multirow{2}*{ \text{PWLS} }  & \text{Error}   &  3.70e-6 & 6.05e-6   & 4.47e-5  & 4.23e-5     \\
& \text{Time}   &  1.34e+1 &   4.03e+1   &  4.29e+2   &  1.13e+3     \\ \hline
   \end{tabular}
     \end{center}

\vskip 0.1in

The results listed in Table \ref{2dinhelmcomtable} indicate that the approximation generated by the DPWNN method is much more accurate than that generated by the original PWLS method. %, and the computational cost for the entire iteration of DPWNN is almost twice that of the PWLS method.
%This is mainly caused by adaptive superpositions with corrected plane wave propagation directions generated by iterative algorithms.

\subsection{An example of three-dimensional Maxwell's equations in inhomogeneous media}
Let the boundary data ${\bf g}$ in the Maxwell system
(\ref{maxeq}) be chosen as the same function ${\bf g}_{ex}$ given in section \ref{3dhomomax}. But we choose $\varepsilon$ as
%In this part, we consider another example with complicated structure of $\varepsilon$. We define it as
 \begin{eqnarray}\label{3dinmax}
\varepsilon =\left\{\begin{array}{ll} 2+2i ~~\quad y < 0.5, z < 0.5, \\
1+i ~~\quad y > 0.5, z < 0.5, \\
\frac{1}{2}+\frac{1}{2}i ~~\quad y > 0.5, z > 0.5, \\
\frac{3}{2}+\frac{3}{2}i ~~\quad y < 0.5, z > 0.5.
\end{array}\right.
\end{eqnarray}
The other quantities are kept the same. %Since the parameter $\varepsilon$ is different from the one chosen in section \ref{3dhomomax}, the vector function ${\bf E}_{\text{ex}}$ in section \ref{3dhomomax} is not the analytic solution of the current equations.

% For the case of inhomogeneous media (i.e., $\omega$ is not a constant), it is difficult to give an analytic solution of the Helmholtz system (\ref{helm1}). As usual we replace the analytic solution by its ``good" approximation generated by the standard finite element method with very fine grids in order to compute accuracies of the plane wave approximations generated by the proposed method.

\begin{figure}[H]
%\vspace{-2cm}thb
\begin{center}
\begin{tabular}{cc}
\epsfxsize=0.5\textwidth\epsffile{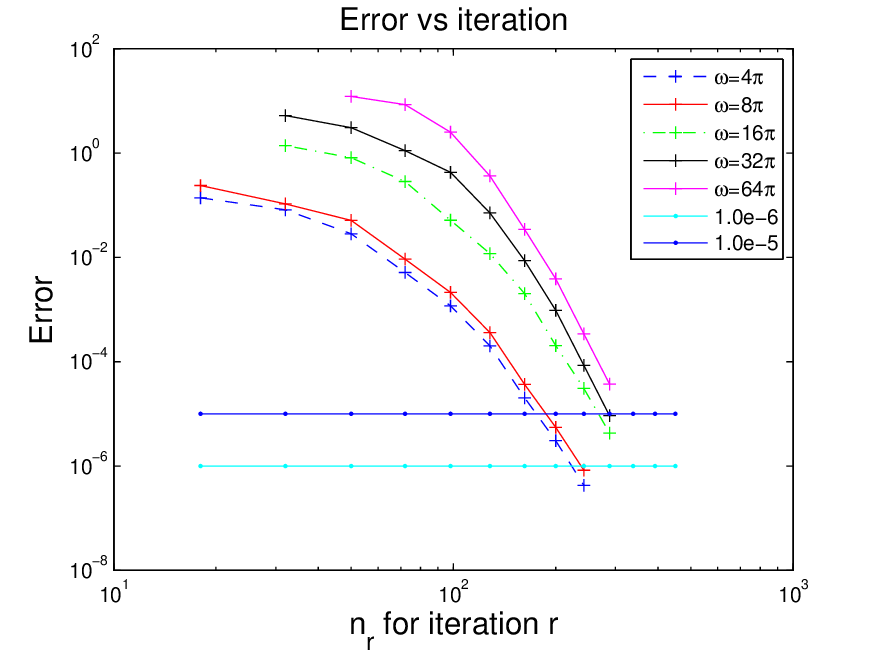}&
\epsfxsize=0.5\textwidth\epsffile{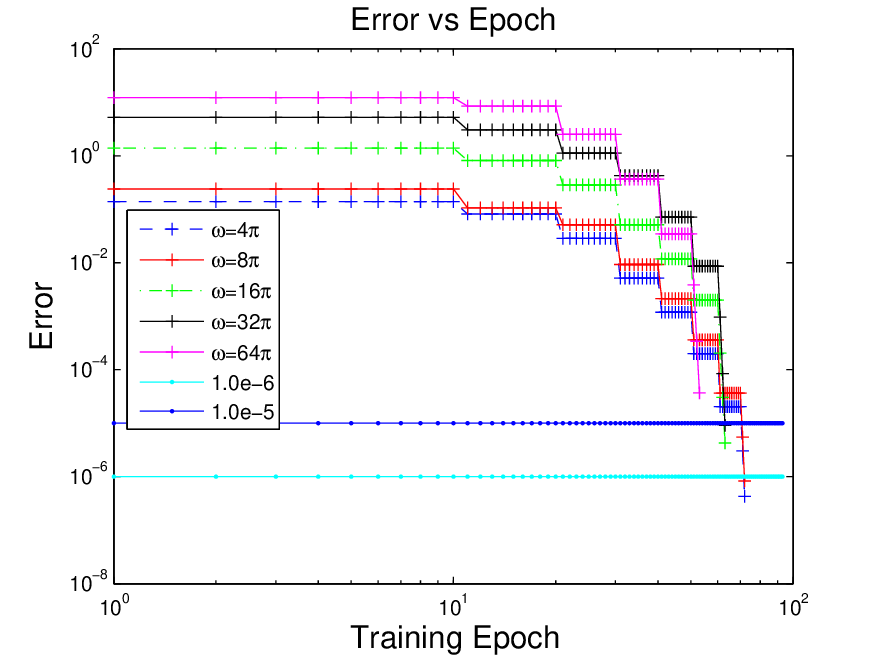}\\
\end{tabular}
\end{center}
 \caption{Displacement of a string (Maxwell's equations in three dimensions). (Left) Errors at each quasi-minimization iteration. (Right) The progress of the loss function within each quasi-minimization iteration.  }
\label{3dinhomo_mul_larwave_nn1}
\end{figure}

Figure \ref{3dinhomo_mul_larwave_nn1} shows the estimated errors of the true errors $|||{\bf E}-{\bf E}_{r-1}|||$. %at the end of each quasi-minimization iteration. We also provide the analogous results after each training epoch.
 The relative parameters are set as follows.
\beq \nonumber
\omega& = & 4\pi, 8\pi, ~ h=\frac{1}{2}, ~m^{\ast}_{r} = r+2, ~tol = 10^{-6}; \cr
\omega& = & 16\pi, ~~~~~ h=\frac{1}{4}, ~m^{\ast}_{r} = r+3, ~tol = 10^{-5};\cr
\omega& = & 32\pi, ~~~~~ h=\frac{1}{8},~m^{\ast}_{r} = r+3, ~tol = 10^{-5};\cr
\omega& = & 64\pi, ~~~~~ h=\frac{1}{16},~m^{\ast}_{r} = r+4, ~tol = 10^{-4}.
\eq
%In general, we choose $h\approx \mathcal{O}(\frac{4\pi}{\omega})$. Here we need to use relatively great $n_r$ since we have chosen relatively coarse meshes.

Numerical results validate that, the resulting approximate solutions
can reach the given accuracy. %when choosing $h\approx \mathcal{O}(\frac{4\pi}{\omega})$ and gradually increasing the width $n_r$ of the network at every iteration step of Algorithm 4.1,

%Next, we would like to compare the proposed DPWNN with the PWLS method. %Let the number of degree of freedoms per each element in the PWLS method be the same as the width of the discontinuous network in the final step of Algorithm 4.1. We set the maximum number of training epoches as 2 in Algorithm 4.2.
  Table \ref{3dinmaxcomtable} shows the comparison of the errors of the resulting approximations and the computing time between the proposed DPWNN and the PWLS method.

\vskip 0.1in
\begin{center}
       \tabcaption{}\vskip -0.3in
\label{3dinmaxcomtable}
       Comparison of errors of approximations and the computing time with respect to $\omega$.  \vskip 0.1in
\begin{tabular}{|c|c|c|c|c|c|c|} \hline
  \multicolumn{2}{|c| } {  \(\omega\) } & $4\pi$ & $8\pi$ &  $16\pi$ & $32\pi$ & $64\pi$ \\ \hline
\multirow{2}*{ \text{DPWNN} } & \text{Error}   &  4.36e-7  & 8.26e-7   &  4.95e-6  & 9.39e-6 & 4.03e-5 \\
& \text{Time}   & 3.03e+2   &  5.69e+2  & 6.14e+3  & 1.62e+4  & 1.36e+5 \\ \hline
 \multirow{2}*{ \text{PWLS} }  & \text{Error}   & 2.47e-6 &   4.87e-6 &  2.99e-5  & 6.45e-5  &  2.93e-4 \\
& \text{Time}   &  1.59e+2  & 2.97e+2   & 3.21e+3  & 7.97e+3  &  6.91e+4  \\ \hline
   \end{tabular}
     \end{center}

\vskip 0.1in

The results listed in Table \ref{3dinmaxcomtable} indicate that the approximation generated by the DPWNN method is much more accurate than that generated by the original PWLS method.
%This is mainly caused by adaptive superpositions with corrected plane wave propagation directions generated by iterative algorithms.

\section{Summary}

We have introduced a {\it discontinuous} plane wave neural network method with $hp-$refinement for approximately solving Helmholtz equation and Maxwell's equations. In this method, we define a  quadratic minimization problem
 and introduce new discretization sets spanned by element-wise neural network functions, where the activation function is chosen as the plane wave-type function. The desired approximate solutions are recursively generated by
 iteratively solving the quasi-minimization problems on the discontinuous plane wave neural networks with a single hidden layer, where plane wave direction angles and activation coefficients are alternatively computed by the proposed iterative algorithms. We emphasize that the performance of the algorithms is mildly dependent on wave numbers, by choosing $h\approx \mathcal{O}(\frac{4\pi}{\omega})$ and gradually increasing the width $n_r$ of the networks.
%, and only no more than ten minimization iterations may guarantee the convergence of the discontinuous plane wave neural network algorithm.

%In the future, we shall focus on the extension of the proposed DPWNN to the models of nonhomogeneous wave equations, anisotropic wave equations and time-dependent wave equations.

\end{document}